\documentclass[a4paper,11pt]{article}
\usepackage{amsmath,amssymb,amsthm,psfrag,color,graphicx,caption2}

\renewcommand{\captionlabeldelim}{~}

\newtheorem*{theoremA}{Theorem A}
\newtheorem*{theoremB}{Theorem B}
\newtheorem*{theoremC}{Theorem C(The straightening theorem)}
\newtheorem*{theoremD}{Theorem D}
\newtheorem*{Theorem}{Theorem}

\newtheorem{lemma}{Lemma}
\newtheorem*{corollary}{Corollary}
\newtheorem{prop}{Proposition}
\newtheorem*{mainth}{Main Theorem}
\newtheorem*{mainprop}{Main Proposition}
\newtheorem*{KLlemma}{Kahn-Lyubich Lemma}
\newtheorem*{claim1}{Claim 1}
\newtheorem*{claim2}{Claim 2}
\newtheorem*{claim3}{Claim 3}
\newtheorem*{subclaim}{Subclaim}

\theoremstyle{definition}
\newtheorem{definition}{Definition}
\newtheorem*{remark}{Remark}

\newtheorem*{acknowledgements}{Acknowledgements}

\numberwithin{equation}{section}

\begin{document}

\title{Proof of the Branner-Hubbard conjecture on Cantor Julia sets}

\author{Weiyuan Qiu \and Yongcheng Yin}
\date{}
\maketitle

\begin{abstract}
By means of a nested sequence  of some critical pieces constructed
by Kozlovski, Shen, and van Strien, and by using a covering lemma
recently proved by Kahn and Lyubich, we prove that the Julia set
of a polynomial is a Cantor set if and only if each component of
the filled-in Julia set containing critical points is aperiodic.
This result was a conjecture raised by Branner and Hubbard in
1992.
\end{abstract}

\section{Introduction}\label{intro}

\noindent For a complex polynomial $f$ of degree $d\geqslant 2$, the
set
$$K_f = \{z\in \mathbb{C}\,\,| \, \mbox{ the sequence } \{f^n(z)\} \mbox{ is bounded}\}$$
is called the filled-in Julia set of $f$, where $f^n$ is the
$n$-th iterate of $f$. The Julia set $J_f$ of $f$ is the boundary
of $K_f$. A component of $K_f$ is called critical if it contains
critical points. We denote the component of $K_{f}$ containing $x$
by $K_{f}(x)$. A component $K_{f}(x)$ is aperiodic if $f^n K_f
(x)\not=K_f(x)$ for all $n>0$.

P.\,Fatou and G.\,Julia proved the following theorem.

\begin{theoremA} [\cite{F} and \cite{J}] (1)\; The Julia set of a
complex polynomial $f$ is connected if and only if $K_f$ contains
all critical points of $f$.

(2)\; The Julia set of a complex polynomial $f$ is a Cantor set if
$K_f$ contains no critical points of $f$.
\end{theoremA}

Fatou conjectured that the condition in Theorem A(2) is also
necessary for the Julia set to be a Cantor set. But this was
disproved by Brolin in \cite{B}. He gave some real cubic polynomials
with Cantor Julia set $J_f=K_f$ containing one critical point.

Using combinatorial system of tableaus, Branner and Hubbard
completely settled the question of when the Julia set of a cubic
polynomial is a Cantor set. They proved

\begin{theoremB}[\cite{BH-2}]
For a cubic polynomial $f$ with one critical point in $K_f$, the
Julia set $J_f$ is a Cantor set if and only if the critical
component of $K_f$ is aperiodic.
\end{theoremB}

The same combinatorics was used by Yoccoz to prove the local
connectivity of the Julia set of a quadratic polynomial which has
no irrational indifferent periodic points and which is not
infinitely renormalizable. Transferring this result to parameter
space, he proved that the Mandelbrot set is locally connected at
these parameters. See \cite{H} and \cite{M}. Yoccoz introduced a
partition of the complex plane by using external rays and
equipotential curves. Such partition is called a {\it{Yoccoz
puzzle}}. It becomes a powerful tool in the study of dynamics of
polynomials, see for example \cite{AKLS}, \cite{GS}, \cite{H},
\cite{Ji}, \cite{KL-2}, \cite{KSS}, \cite{LS}, \cite{Ly-2},
\cite{LY}, \cite{Mc-1}, \cite{Mc-2}, \cite{M}, \cite{Sh}, and
\cite{Yo}. In \cite{Ji}, Jiang gives the first proof that the
Julia set of an unbranched infinitely renormalizable quadratic
polynomial having complex bounds is locally connected. A different
proof has been given by McMullen in \cite{Mc-1}. Other puzzles are
used to prove local connectivity of the Julia sets of some
quadratic Siegel polynomials and cubic Newton maps, see \cite{Pe},
\cite{PZ}, \cite{Ro-1} and \cite{Ro-2}.

In \cite{BH-2}, Branner and Hubbard conjectured that the assertion
in theorem B is true for any polynomial.

Let $f$ be a polynomial with real coefficients such that one real
critical point has a bounded orbit and all other critical points
escape to infinity. Then the Julia set $J_f$ is a Cantor set if
and only if  the critical component of $K_f$ is aperiodic. See
\cite{LS1} and \cite{LS2}. In \cite{E}, Emerson gave a
combinatorial condition for the Julia set of a polynomial to be a
Cantor set and showed that there are polynomials fulfilling the
condition.

The purpose of this paper is to give a proof of the above
Branner-Hubbard's conjecture. We state the main result of this
paper as the

\begin{mainth} Let $f$ be a complex polynomial of degree $\geqslant 2$ and
let $\mathrm{Crit}$ be the set of critical points of $f$ with
bounded orbits. Then the Julia set $J_f$ of $f$ is a Cantor set if
and only if the critical component $K_f(c)$ is aperiodic for all
$c\in$$\mathrm{Crit}$.
\end{mainth}

There are two important tools in our proof. One is a nested
sequence of some critical pieces constructed by Kozlovski, Shen,
and van Strien in \cite{KSS} which we shall call ``KSS nest". The
other one is a covering lemma proved by Kahn and Lyubich recently,
see \cite{KL-1}. This covering lemma has many important
applications in complex dynamics, see \cite{AKLS} and \cite{KL-2}.

The paper is organized as follows. In section 2, we present some
definitions and reduce the Main Theorem to the Main Proposition.
We summarize the construction of KSS nest in section 3. The proof
of the Main Proposition is given in section 4. In section 5, we
prove a stronger result than the Main Theorem which states that
each wandering component of the filled-in Julia set for an
arbitrary polynomial is a point.

\section{Definitions and preliminary results}\label{defi_resu}

\noindent For a complex polynomial $f$ of degree $d\geqslant 2$,
it is well-known that the function
$$G:\mathbb{C}\rightarrow \mathbb{R}_{+}\cup\{0\}$$
defined by
$$G(z)=\lim_{n\to\infty} \frac{1}{d^n} \log ^+ |f^n(z)|$$
is continuous and satisfies
$$(1)\quad G(f(z))=d G(z), \qquad (2)\quad K_f=\{z\in \mathbb{C}\,|\, \,G(z)=0\},$$
see \cite{BH-1} and \cite{DH-1}.

The {\it{Branner-Hubbard puzzle}} of $f$ is constructed as
follows. Choose a small number $r_0>0$ which is not a critical
value of $G$ such that the region $G^{-1}(0,r_0)$ contains no
critical points of $f$. Then for each integer $k\geqslant 0$, the
locus
$$G^{-1}([0,r_0 d^{-k})) = \{z\in\mathbb{C} |\, G(z)< r_0 d^{-k}\}$$
is the disjoint union of a finite number of open topological
disks. Each such open disk will be called a puzzle piece $P_k$ of
depth $k$. Thus each point $x\in K_f$ determines a nested sequence
$P_0(x)\supset P_1(x)\supset \cdots$ and
$K_f(x)=\bigcap_{k\geqslant 0}P_k(x)$. By Gr\"{o}tzsch's
inequality,
$$\mod(P_{0}(x) - K_{f}(x)) = \infty$$
if there exists a subsequence
$$P_{k_1^\prime}(x)\supset P_{k_1}(x)\supset
P_{k_2^\prime}(x)\supset P_{k_2}(x)\supset\cdots$$ such that
$$\sum_{1=1}^\infty \!\!\!\mod(P_{k_i^\prime}(x) - \overline{P_{k_i}(x)}) = \infty.$$
It follows that $K_f(x)=\bigcap_{k\geqslant 0}P_k(x)= \{x\}$, see
\cite{Ah} and \cite{BH-2}.

The Julia set $J_f$ is a Cantor set if and only if $K_f(x)=\{x\}$
for any $x\in K_f$.

If a component $K$ of $K_f$ contains critical points
$c_1,c_2,\cdots,c_k$, then $P_n(c_1)=P_n(c_2)=\cdots =P_n(c_k)$
and
$$\deg (f|_{P_n(c_1)})=(\deg_{c_1}f-1)+(\deg_{c_2}f-1)+\cdots+(\deg_{c_k}f-1)+1$$
for all $n\geqslant 0$. We can think of $K$ as a component
containing one critical point of degree
$(\deg_{c_1}f-1)+(\deg_{c_2}f-1)+\cdots+(\deg_{c_k}f-1)+1$. We
therefore assume each critical component of $K_f$ contains only
one critical point in the following. Take $r_0$ small enough such
that each puzzle piece contains at most one critical point.

For each $x\in K_f$, the tableau $T(x)$ is defined in \cite{BH-2}.
It is the two dimension array $P_{n,l} (x) = f^l (P_{n+l}(x))$.
The position $(n,l)$ is called critical if $P_{n,l}(x)$ contains a
critical point of $f$. If $P_{n,l}(x)$ contains a critical point
$c$, the position $(n,l)$ is called a $c$-position. Let
$\mathrm{Crit}$ be the set of critical points with bounded orbits.
The tableau $T(c)$ of a critical point $c\in$ $\mathrm{Crit}$ is
called periodic if there is a positive integer $k$ such that
$P_n(c) = f^k (P_{n+k}(c))$ for all $n\geq 0$. Otherwise, $T(c)$
is said to be aperiodic.

All the tableaus satisfy the following three rules
\begin{enumerate}
\item[(T1)] If $P_{n,l}(x)=P_n(c)$ for some critical point $c$,
then $P_{i,l}(x)=P_i(c)$ for all $0\leqslant i \leqslant n$.

\item[(T2)] If $P_{n,l}(x)=P_n(c)$ for some critical point $c$,
then $P_{i,l+j}(x)=P_{i,j}(c)$ for $i+j\leqslant n$.

\item[(T3)] Let $T(c)$ be a tableau for some critical point $c$
and $T(x)$ be any tableau. Assume
           \begin{enumerate}
           \item $P_{n+1-l,l}(c)=P_{n+1-l}(c_1)$ for some critical
                 point $c_1$ and $n>l\geqslant
                 0$, and $P_{n-i,i}(c)$ contains no critical points for
                 $0<i<l$.
           \item $P_{n,m}(x)=P_n(c)$ and $P_{n+1,m}(x)\neq
                  P_{n+1}(c)$ for some $m>0$.
           \end{enumerate}
Then $P_{n+1-l,m+l}(x)\neq P_{n+1-l}(c_1)$.
\end{enumerate}

In order to show that the Julia set for a polynomial is a Cantor
set, we shall use the polynomial-like mapping theory introduced by
Douady and Hubbard in \cite{DH-2}. Recall that a
{\it{polynomial-like mapping}} of degree $d$ is a triple $(U,V,g)$
where $U$ and $V$ are simply connected plane domains with
$\overline{V}\subset U$, and $g:V\rightarrow U$ is a holomorphic
proper mapping of degree $d$. The filled-in Julia set $K_g$ of the
polynomial-like mapping $g$ is defined as
$$K_g=\{z\in V\,|\,\,g^n(z)\in U \,\,\mathrm{for \,\,all}\,\, n\geqslant0\}.$$

Two polynomial-like mappings $(U_1,V_1,g_1)$ and $(U_2,V_2,g_2)$
of degree $d$ are said to be {\it{hybrid equivalent}} if there
exists a quasi-conformal homeomorphsim $h$ from a neighborhood of
$K_{g_1}$ onto a neighborhood of $K_{g_2}$, conjugating $g_1$ and
$g_2$ and such that $\bar{{\partial}}h=0$ on $K_{g_1}$. The
following theorem was proved by Douady and Hubbard in \cite{DH-2}.

\begin{theoremC}  (1)\; Every polynomial-like
mapping $(U,V,g)$ of degree $d$ is hybrid equivalent to a
polynomial of degree $d$.

(2)\; If $K_g$ is connected, then the polynomial is uniquely
determined up to conjugation by an affine map.
\end{theoremC}

If $T(c)$ is periodic of period $k$, then
$(P_n(c),P_{n+k}(c),f^k)$ is a polynomial-like mapping of degree
$\deg(f^k|_{K_f(c)})\geqslant 2$ for some $n\geqslant 0$. The
filled-in Julia set of this polynomial-like mapping equals to
$K_f(c)$. From the straightening theorem, $K_f(c)$ is
quasi-conformally homeomorphic to the filled-in Julia set of a
polynomial of degree $\deg(f^k|_{K_f(c)})\geqslant 2$. The ``only
if" part in the Main Theorem is obvious. We always assume that
each critical component of $K_f$ is aperiodic before section 5. It
is equivalent to assuming $T(c)$ is aperiodic for all $c\in
\mathrm{Crit}$.

\begin{definition}
(1)\; The tableau $T(x)$ for $x\in K_f$ is {\it{non-critical}} if
there exists an integer $n_0\geqslant 0$ such that $(n_0,j)$ is not
critical for all $j>0$.

(2)\; We say the forward orbit of $x$ {\it{combinatorially
accumulates}} to $y$, written as $x\to y$, if for any $n\geqslant
0$, there exists $j>0$ such that $y\in P_{n,j}(x)$, i.e.,
$f^j(P_{n+j}(x))= P_n(y)$. It is clear that if $x\to y$ and $y\to
z$, then $x\to z$. For each critical point $c\in\mathrm{Crit}$,
let
$$F(c)=\{c^\prime\in\mathrm{Crit}\,| \,\, c\to c^\prime\}$$
and
$$[c]=\{c^\prime\in\mathrm{Crit}\,| \,\, c\to c^\prime \textrm{ and } c^\prime \to
c \}.$$

(3)\;We say $P_{n+k}(c^\prime)$ is a {\it{child}} of $P_n(c)$ if
$c^\prime\in [c]$, $f^k(P_{n+k}(c^\prime))= P_n(c)$, and $f^{k-1}:\,
P_{n+k-1}(f(c^\prime))\to P_n(c)$ is conformal.

(4)\; Suppose $c\to c$ , i.e., $[c]\not=\emptyset$. We say $T(c)$
is {\it{persistently recurrent}} if $P_n(c_1)$ has only finitely
many children for all $n\geqslant 0$ and all $c_1\in [c]$.
Otherwise, $T(c)$ is said to be {\it{reluctantly recurrent}}.
\end{definition}

Take a small $r_0>0$ such that for any $c\in\mathrm{Crit}$, there
are no $c^\prime$-positions in the first row of $T(c)$ if
$c\not\to c^\prime$.

\begin{figure}[h]
\psfrag{c}{$c$}\psfrag{c1}{$c_1$}\psfrag{c2}{$c_2$} \centering
\includegraphics[width=10cm]{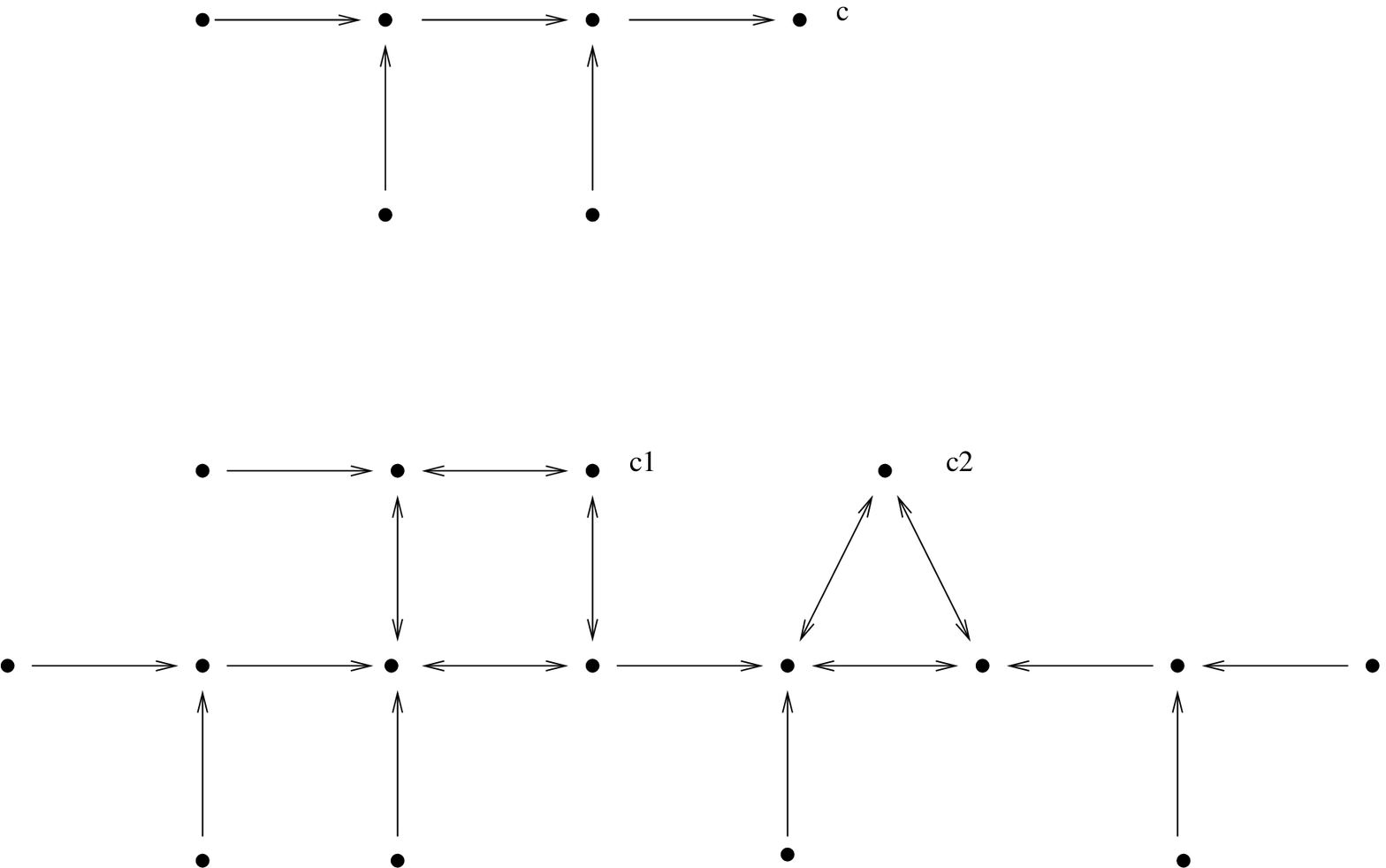}\caption{}\label{fig1}
\end{figure}

Let
\begin{eqnarray*}
&&\textrm{Crit}_\textrm{n}=\{c\in\mathrm{Crit}\,| \,\, T(c) \textrm { is non-critical}\},\\
&&\textrm{Crit}_\textrm{p}=\{c\in\mathrm{Crit}\,| \,\, T(c) \textrm { is persistently recurrent}\},\\
&&\textrm{Crit}_\textrm{r}=\{c\in\mathrm{Crit}\,| \,\, T(c) \textrm { is reluctantly recurrent}\},\\
&&\textrm{Crit}_\textrm{en}=\{c^\prime\in\mathrm{Crit}\,|\, \,
c^\prime \not\to c^\prime
\textrm{ and } c^\prime \to c \textrm{ for some } c\in\textrm{Crit}_\textrm{n} \},\\
&&\textrm{Crit}_\textrm{ep}=\{c^\prime\in\mathrm{Crit}\,| \,\,
c^\prime \not\to c^\prime
\textrm{ and } c^\prime \to c \textrm{ for some } c\in \textrm{Crit}_\textrm{p}\},\\
&&\textrm{Crit}_\textrm{er}=\{c^\prime\in\mathrm{Crit}\,| \,\,
c^\prime \not\to c^\prime \textrm{ and } c^\prime \to c \textrm{
for some } c\in \textrm{Crit}_\textrm{r}\}.
\end{eqnarray*}
Then
$$\mathrm{Crit}=\textrm{Crit}_\textrm{n}\cup \textrm{Crit}_\textrm{p}\cup \textrm{Crit}_\textrm{r}\cup \textrm{Crit}_\textrm{en}
\cup \textrm{Crit}_\textrm{ep}\cup \textrm{Crit}_\textrm{er}.$$ It
is not a classification because these sets might intersect.

Consider the critical points $c$, $c_1$ and $c_2$ in
Figure~\ref{fig1}. The tableau $T(c)$ for $c$ is non-critical. From
Lemma 1 in the following, the tableau $T(c_{1})$ for $c_1$ is
reluctantly recurrent. The tableau $T(c_{2})$ for $c_2$ is
reluctantly recurrent or persistently recurrent.

Combined with arguments of Branner and Hubbard in \cite{BH-2}, we
have the following proposition.

\begin{prop}
(1)\; If $T(x)$ is non-critical, then $K_f(x)=\{x\}$.

(2)\; Suppose $c\in$$\mathrm{Crit}$$_\mathrm{n}\cup$
$\mathrm{Crit}$$_\mathrm{r}$. Then $K_f(c)=\{c\}$ and
$K_f(x)=\{x\}$ for any $x\to c$.

\end{prop}

\begin{proof}
Let $\mathcal{P}_n$ be the collection of all puzzle pieces of
depth $n$. It has only finitely many pieces. Hence
$$\nu_n=\min\{\mathrm{mod}(P_n - \overline{P_{n+1}})| \,P_n \in \mathcal{P}_n,\,
P_{n+1}\in \mathcal{P}_{n+1}, \mathrm{with}\,\,P_{n+1}\subset
P_n\}>\,0.$$

(1) Since $T(x)$ is non-critical, there exists an integer
$n_0\geqslant 0$ such that $(n_0,j)$ is not a critical position for
all $j>0$. For any $k\geqslant1$,
$\deg(f^k|_{P_{n_0+k}(x)})\leqslant \deg(f|_{P_{n_0+1}(x)})$ and
\begin{equation}
\begin{split}
\mathrm{mod}(P_{n_0+k}(x)- \overline{P_{n_0+k+1}(x)}) &\geqslant \frac{1}{\deg(f^k|_{P_{n_0+k}(x)})}\mathrm{mod}(P_{n_0}(f^k(x)) - \overline{P_{n_0+1}(f^k(x))})\\
                                   &\geqslant \frac{\nu_{n_0}}{\deg(f|_{P_{n_0+1}(x)})}>0.
\end{split}
\end{equation}
This yields
$$\sum_{n=1}^\infty \!\!\!\mod(P_n(x) - \overline{P_{n+1}(x)}) =
\infty.$$
Hence, $K_f(x)=\{x\}$.

(2) If $c \in \mathrm{Crit}_\mathrm{n}$, then $T(c)$ is non-critical
and $K_f(c)=\{c\}.$ There exists an integer $n_0\geq 0$ such that
$$\deg(f^k|_{P_{n_0+k}(c)})=\deg_cf$$
for all $k\geqslant 1$.

For any $x\to c$, let $l_k$ be the first moment such that
$f^{l_k}(x)\in P_{n_0+k}(c)$, $i.e.$ $(n_0+k, l_k)$ is the first
$c$-position on the $(n_0+k)$-th row in $T(x)$. By tableau rules
(T1) and (T2), there is at most one $c^\prime$-position on the
diagonal
$$\{(n,m)\,|\, \, n+m = n_0+k+l_k, \quad n_0+k < n \leqslant n_0+k+l_k \}$$
for any $c^\prime \in \mathrm{Crit} - \{c\}.$ Therefore,
$f^{l_k+k}(P_{n_0+k+l_k}(x))= P_{n_0}(f^{l_k+k}(x))$ and $\deg
(f^{l_k+k}|_{P_{n_0+k+l_k}(x)})\leqslant D_1 < \infty$ for any
$k\geq 1$, where $D_1$ is an integer independent of $k$. We have
$$\mod(P_{n_0+k+l_k}(x) - \overline{P_{n_0+k+l_k+1}(x)})\geqslant \frac{\nu_{n_0}}{D_1}>0.$$
So $K_f(x)=\{x\}$.

If $c \in \mathrm{Crit}_\mathrm{r}$, then there exist an integer
$n_0\geq 0$, $c^\prime \in [c]$,  $c_1 \in [c]$ and infinitely
many integers $k_n\geqslant 1$ such that
$\{P_{n_0+k_n}(c^\prime)\}_{n\geqslant 1}$ are children of
$P_{n_0}(c_1)$. Let $m_n$ be the first moment such that
$f^{m_n}(c)\in P_{n_0+k_n}(c^\prime)$. There is at most one
$\widetilde{c}$-position on the diagonal
$$\{(n,m)\,|\,\, n+m = n_0+k_n+m_n, \quad n_0+k_n < n \leqslant n_0+k_n+m_n\}$$
in $T(c)$ for any $\widetilde{c}\in \mathrm{Crit} - \{c\}.$
Therefore, $f^{m_n+k_n}(P_{n_0+k_n+m_n}(c))= P_{n_0}(c_1)$ and
$\deg (f^{m_n+k_n}|_{P_{n_0+k_n+m_n}(c)})\leqslant D_2 < \infty$
for any $n\geq 1$, where $D_2$ is an integer independent of $n$.
We have
$$\mod(P_{n_0+k_n+m_n}(c) - \overline{P_{n_0+k_n+m_n+1}(c)})\geqslant \frac{\nu_{n_0}}{D_2}>0$$
and $K_f(c)=\{c\}$.

Suppose $x\to c$ for some $c \in \mathrm{Crit}_\mathrm{r}$. Let
$l_n$ be the first moment such that $f^{l_n}(x)\in
P_{n_0+k_n+m_n}(c)$ and let $t_n = k_n+m_n+l_n$. By the same
method, we have $f^{t_n}(P_{n_0+t_n}(x))= P_{n_0}(c_1)$ and $\deg
(f^{t_n}|_{P_{n_0+t_n}(x)})\leqslant D_3 < \infty$ for any $n\geq
1$, where $D_3$ is an integer independent of $n$. Hence
$$\mod(P_{n_0+t_n}(x) - \overline{P_{n_0+t_n+1}(x)})\geqslant \frac{\nu_{n_0}}{D_3}>0$$
and $K_f(x)=\{x\}$.
\end{proof}

From Proposition 1 and $\mathrm{Crit}=\textrm{Crit}_\textrm{n}\cup
\textrm{Crit}_\textrm{p}\cup \textrm{Crit}_\textrm{r}\cup
\textrm{Crit}_\textrm{en} \cup \textrm{Crit}_\textrm{ep}\cup
\textrm{Crit}_\textrm{er},$ we can reduce the Main Theorem to the
following proposition.

\begin{mainprop}
If $c\in $$\mathrm{Crit}$$_\mathrm{p}$, then $K_f(c)=\{c\}$ and
$K_f(x)=\{x\}$ for all $x\to c$.
\end{mainprop}

The following lemma will be used in sections 3 and 4.

\begin{lemma}\label{lemma1}
If $T(c)$ is persistently recurrent, then  $F(c)=[c]$.
\end{lemma}

\begin{figure}[h]
\psfrag{Tc}{$T(c):$} \psfrag{c}{$c$}\psfrag{n0}{$0$}

\psfrag{c1}{$c_1(k)$}\psfrag{cp}{$c'$}\psfrag{c2}{$c_2(k)$}
\psfrag{0sk}{$(0,s_k)$}\psfrag{tksksk}{$(t_k-s_k,s_k)$}
\psfrag{n0mik}{$(0,m_{k})$} 
\psfrag{0nkmk}{$(0,n_k+m_k)$}\psfrag{nkmk}{$(n_k,m_k)$}
\psfrag{mik}{$(m_{i_k+1}-m_{i_k},m_{i_k})$}
\psfrag{0tk}{$(0,t_k)$}
\centering\includegraphics[width=10cm]{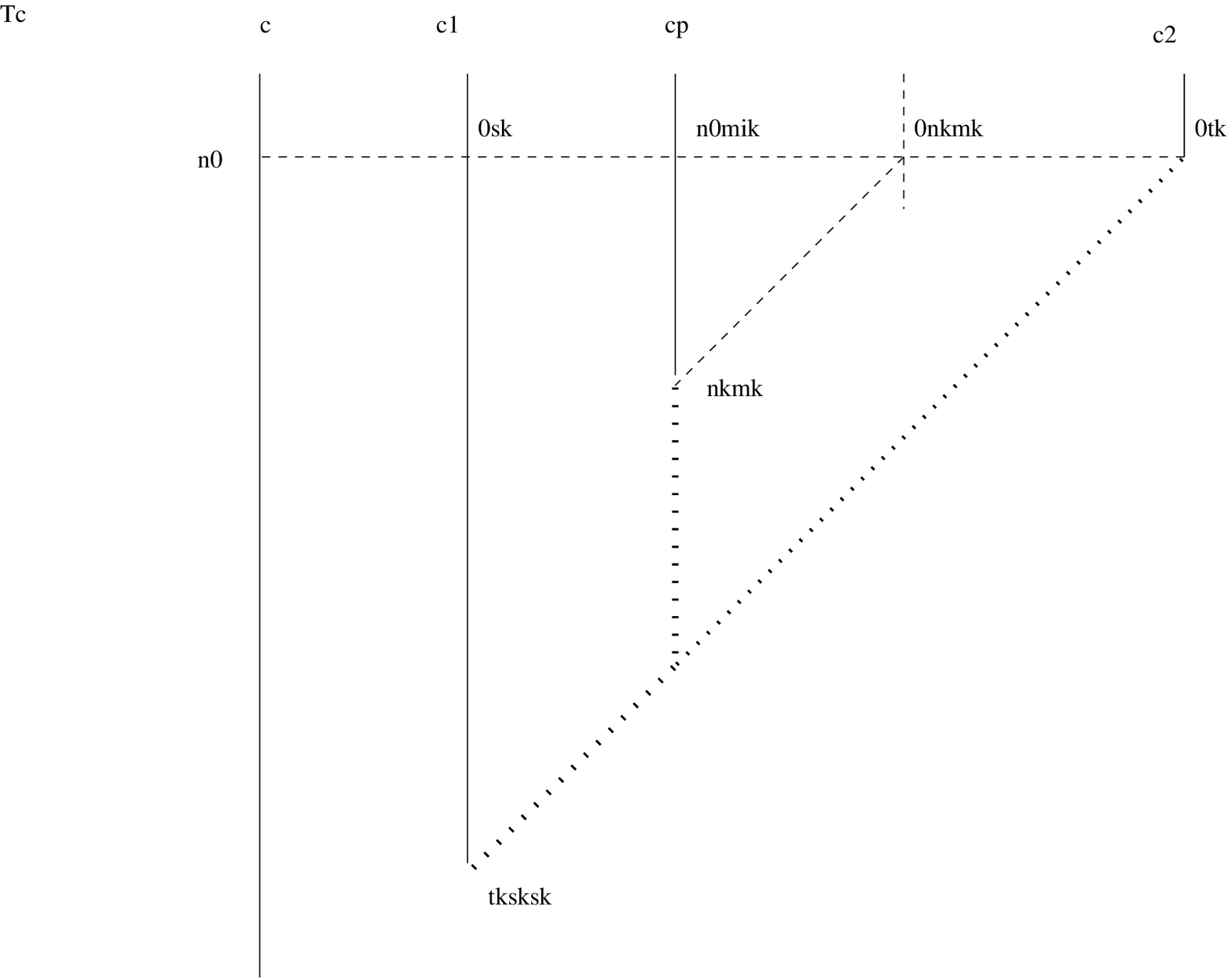}\caption{}\label{fig2}
\end{figure}

\begin{proof}
Suppose $c\to c^\prime$ and $c^\prime\not\to c$. If there exists a
column where each position is $c'$-position, then $c^\prime \to
c$. It contradicts with our assumption. Hence there are infinitely
many $c^\prime$-positions $\{(n_k,m_k)\}_{k\geqslant 1}$ in $T(c)$
such that $(n_k+1,m_k)$ is not critical and $\lim_{k\to \infty}n_k
= \infty$. By the tableau rule (T2) and the choice of $r_0$, there
are no $\tilde{c}$-positions on the diagonal
$$\{(n,m)\,|\,\, n+m = n_k+m_k, \quad 0\leqslant n \leqslant n_k \}$$
for any $\tilde{c}\in [c]$.

Let $(0,t_k)$ be a $c_2(k)$-position on the right of $(0,n_k+m_k)$
for some $c_2(k)\in [c]$ such that there are no
$\tilde{c}$-positions between $(0,n_k+m_k)$ and $(0,t_k)$ for any
$\tilde{c}\in [c]$. Then there are no $c^{\prime\prime}$-positions
on the diagonal
$$\{(n,m)\,|\,\, n+m = t_k, \quad 0< n < t_k-m_k \}$$
for any $c^{\prime\prime}\in F(c)$. Hence all positions on this
diagonal are not critical. Let $s_k$ be the largest integer
between $0$ and $m_k$ such that $(t_k-s_k,s_k)$ is a critical
position. Say it is a $c_1(k)$-position for some $c_1(k)\in [c]$.
See Figure~\ref{fig2}. Take a subsequence $\{k_j\}$ such that
$c_2(k_j)=c_2$ for some $c_2\in [c]$. Then the critical piece
$P_0(c_2)$ has infinitely many children. This is impossible
because $T(c)$ is persistently recurrent.
\end{proof}

The proof of the Main Proposition will be given in section 4.

\section{KSS nest}

\noindent For completeness, we summarize the construction of a
critical nest and related results given by Kozlovski, Shen, and
van Strien in \cite{KSS}. Such nest will be called {\it{KSS
nest}}. Principal nest and modified principal nest are used to
study the dynamics of unicritical polynomials, see \cite{AKLS},
\cite{BSS}, \cite{KL-2}, \cite{Ly-1} and \cite{Ly-2}. In
\cite{Ly-2}, Lyubich proved the linear growth of its ``principal
moduli" for quadratic polynomials. This yields the density of
hyperbolic maps in the real quadratic family. The same result is
also obtained by Graceyk and \'{S}wi\c{a}tek in \cite{GS}. See
also \cite{Mc-1} and \cite{Sh}. Recently, the local connectivity
of Julia sets and combinatorial rigidity for unicritical
polynomials are proved in \cite{KL-2} and \cite{AKLS} by means of
principal nest and modified principal nest. For dynamics of
multimodal maps, see \cite{Shen} and \cite{Sm-1}.

Let $A$ be an open  set and $x\in A$. We denote the connected
component of $A$ containing $x$ by $\textrm{Comp}_x(A)$. Given an
open set $X$ consisting of finitely many puzzle pieces(not
necessarily the same depth) such that $f^n(z)\not\in \overline{X}$
for any $z\in \partial X$ and any $n\geqslant 1$, let
$$D(X)=\{z\in \mathbb{C}\,\,|\,\,\exists k\geqslant 1 \,\,s.t.\, \,f^k(z)\in
X\}.$$ The first entry map
$$R_X:D(X)\rightarrow X$$
is defined by $z\mapsto f^{k(z)}(z)$, where $k(z)\geqslant 1$ is
the smallest integer with $f^{k(z)}(z)\in X$. Let $I$ be a
component of $D(X)$. Then there exists an integer $k$ such that
$k(z)=k$ for any $z\in I$ and $f^k(I)$ is a connected component of
$X$. The orbit
$$\{I,\,f(I),\,\cdots,\,f^{k-1}(I)\}$$
meets each critical point at most once and the degree of $f^k$ on
$I$ is uniformly bounded. For any $z\in D(X)$, let $L_z(X)$ be the
connected component of $D(X)$ containing $z$. We further define
$\hat{L}_z(X)=\textrm{Comp}_z(X)$ for $z\in X$ and
$\hat{L}_z(X)=L_z(X)$ for $z\in D(X)-X$.

Suppose $T(c_0)$ is persistently recurrent, then $F(c_0)=[c_0]$.
Let $$b=\#[c_0],\,\,d_0=\deg_{c_0}f,\,\,d_{max}=\max\{\deg_cf
|\,c\in [c_0]\}$$ and
$$\textrm{orb}([c_0])=\cup_{n\geqslant 0}f^n([c_0]).$$

For any puzzle piece $I$ containing $c_0$, we construct puzzle
pieces $P_c^\prime \subset\subset P_c$ for any $c\in [c_0]$ as
follows. Let $T_0=I$ and $J_0=L_{c_0}(I)$. If $R_I(c^\prime)\in
J_0$ for any $c^\prime \in [c_0]-\{c_0\}$, we take
$P_c=\hat{L}_c(T_0)$ and $P_c^\prime =\hat{L}_c(J_0)$ for any
$c\in [c_0]$. If $R_I(c_1)\not\in J_0$ for some $c_1\in
[c_0]-\{c_0\}$, let $T_1=J_0\cup
\textrm{Comp}_{c_1}(R^{-1}_I(L_{R_I(c_1)}(I)))$ and
$J_1=L_{c_0}(T_1)\cup L_{c_1}(T_1)$. If $R_{T_1}(c^\prime)\in J_1$
for any $c^\prime \in [c_0]-\{c_0,\, c_1\}$, we take
$P_c=\hat{L}_c(T_1)$ and $P_c^\prime =\hat{L}_c(J_1)$ for any
$c\in [c_0]$. Repeating this process, we have $T_m=J_{m-1}\cup
\textrm{Comp}_{c_m}(R^{-1}_{T_{m-1}}(L_{R_{T_{m-1}}(c_m)}(T_{m-1})))$
and $J_m=\bigcup_{0\leqslant i \leqslant m}L_{c_i}(T_m)$ for some
$m<b$ such that $R_{T_m}(c^\prime)\in J_m$ for any $c^\prime \in
[c_0]-\{c_0,\, c_1,\, ...,\,c_m\}$. Let $P_c=\hat{L}_c(T_m)$ and
$P_c^\prime =\hat{L}_c(J_m)$ for any $c\in [c_0]$. These two
pieces $P_c^\prime \subset\subset P_c$ satisfy the following two
properties
\begin{enumerate}
\item[(P1)] There exists an integer $l_c$ such that $f^{l_c}(P_c)= I$, $\deg(f^{l_c}:\; P_c\to I) \leqslant
{d_{max}}^{b^2-b}$ and $\#\{i|\, c_0\in f^i (P_c)$, $0\leqslant
i<l_c\} \leqslant b-1$. The piece $P_c^\prime$ is also a pull-back
of $I$.

\item[(P2)] For each $x\in (P_c-P_c^\prime)\cap \textrm{orb}([c_0])$, there
exist a positive integer $k$, a puzzle piece $V(x)$ containing $x$
and $\widetilde{c}\in [c_0]$ such that $f^k:\; V(x)\to
P_{\widetilde{c}}$ is conformal. In fact, let $k\geqslant 1$ be
the first moment such that $f^k(x)\in \cup_{\widetilde{c}\in
[c_0]}P_{\widetilde{c}}$, $f^k(x)\in P_{\widetilde{c}}$ for some
$\widetilde{c}\in [c_0]$, and let $V(x)$ be the component of
$f^{-k}(P_{\widetilde{c}})$ containing $x$, then $V(x)\subset
P_c-P_c^\prime$ and $f^k:\; V(x)\to P_{\widetilde{c}}$ is
conformal.
\end{enumerate}

Since $T(c_0)$ is persistently recurrent, each $P_c$ has only
finitely many children. Let $Q_c$ be the last child of $P_c$. Then
there exists an integer $v_c\geqslant 1$, largest among all the
children of $P_c$, such that $f^{v_c}(Q_c)= P_c$. The set $Q_c$
contains a critical point $c^\prime \in [c_0]$. Let
$$v=\max\{v_c|\,c\in[c_0]\}.$$
Suppose $v=v_{c_1}$ for some $c_1\in [c_0]$. By (P2) as above and
the maximality of $v$, we have $f^{v} (c^\prime) \in
P_{c_1}^\prime$ and
$$ (Q_{c_1}-Q_{c_1}^\prime) \cap \textrm{orb}([c_0]) = \emptyset,$$
where $Q_{c_1}^\prime$ is the connected component of
$f^{-v}(P_{c_1}^\prime)$ containing $c^\prime$.

\begin{figure}[h]
\psfrag{c}{$c_1$}\psfrag{pc1p}{$P_{c_1}'$}\psfrag{pc1}{$P_{c_1}$}
\psfrag{flc1}{$f^{l_{c_1}}$} \psfrag{c0}{$c_0$}
\psfrag{fc0t}{$f^t(c_0)$} \psfrag{W}{$W$} \psfrag{I}{$I$}
\psfrag{fv}{$f^v$}\psfrag{ft}{$f^t$} \psfrag{fk}{$f^k$}
\psfrag{cp}{$c'$}\psfrag{qc1p}{$Q_{c_1}'$}\psfrag{qc1}{$Q_{c_1}$}
\psfrag{AI}{$A(I)$}\psfrag{BI}{$B(I)$} \centering
\includegraphics[width=8cm]{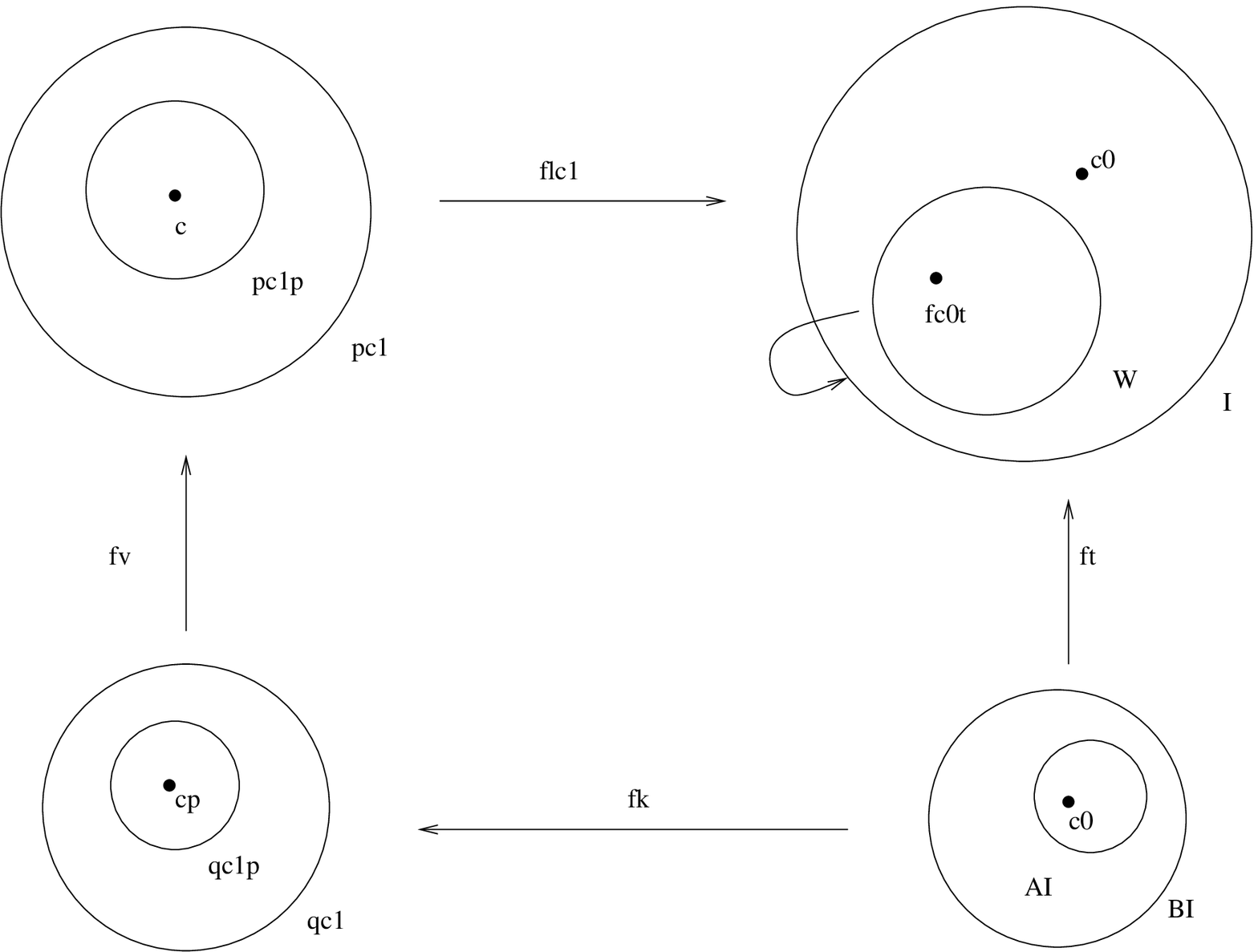}\caption{}\label{fig3}
\end{figure}

\begin{figure}[h]
\psfrag{Tc0}{$T(c_0):$} \psfrag{c0}{$c_0$} \psfrag{I}{$I$}
\psfrag{fs}{$f^s$} \psfrag{ft}{$f^t$} \psfrag{AI}{$A(I)$}
\psfrag{BI}{$B(I)$} \centering
\includegraphics[width=7.5cm]{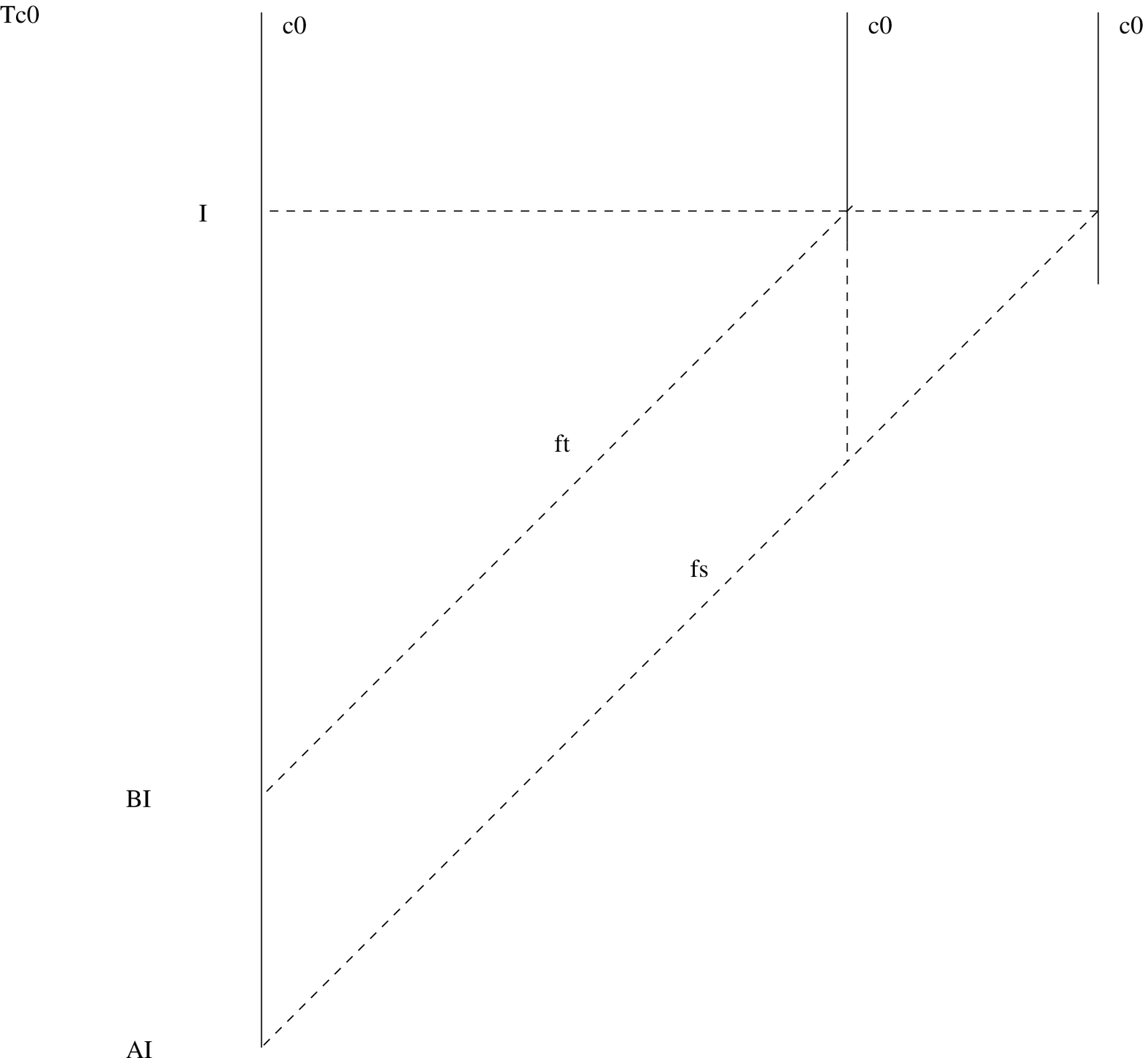}\caption{}\label{fig4}
\end{figure}

Let $B(I)=\hat{L}_{c_0}(Q_{c_1})$, $f^k(B(I))=Q_{c_1}$,
$t=k+v+l_{c_1}$, and $W=L_{f^t(c_0)}(I)$. Then
$f^{l_{c_1}}(P_{c_1}^\prime) \subset W$ because $P_{c_1}^\prime$
is mapped to $I$ and
$$f^t(c_0)=f^{l_{c_1}}(f^{k+v}(c_0))\in f^{l_{c_1}}(P_{c_1}^\prime) \cap W
\not=\emptyset.$$ Let $A(I)$ be the connected component of
$f^{-t}(W)$ containing $c_0$ and $f^s(A(I))= I$. See
Figures~\ref{fig3} and \ref{fig4}.

\begin{figure}[h]
\psfrag{Lcp}{$\hat{L}_c(P)$}\psfrag{c0}{$c_0$}

\psfrag{P}{$P$} \psfrag{Child}{Child}\psfrag{c}{$c$}

\psfrag{fq}{$f^q$}
\psfrag{Lc0q}{$\hat{L}_{c_0}(Q)$}\psfrag{cp}{$c'$}

\psfrag{Q}{$Q$} \centering\includegraphics[width=10cm]{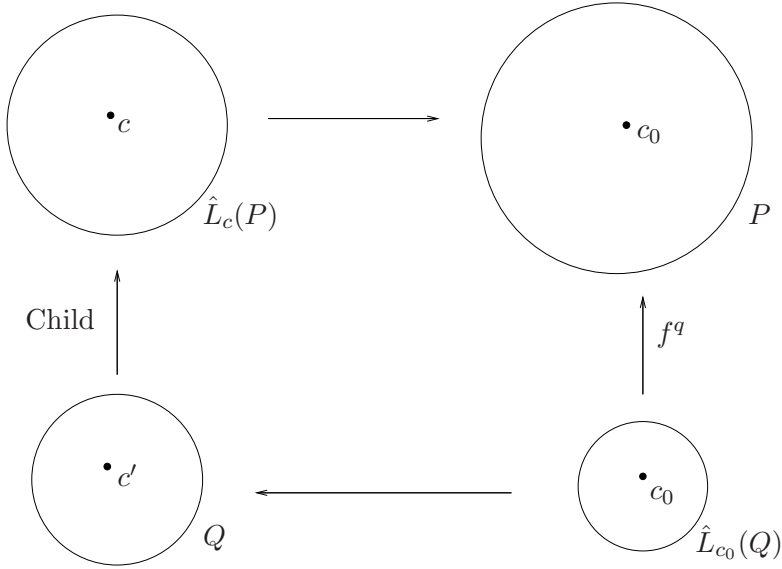}
\caption{A successor of $P$}\label{fig5}
\end{figure}

\begin{definition}
Given  a puzzle piece $P$ containing $c_0$, a {\it{successor}} of
$P$ is a piece of the form $\hat{L}_{c_0}(Q)$, where $Q$ is a
child of $\hat{L}_c(P)$ for some $c\in [c_0]$. See
Figure~\ref{fig5}.
\end{definition}

\renewcommand{\captionlabeldelim}{.}

It is clear that $L_{c_0}(P)$ is a successor of $P$. Since
$T(c_0)$ is aperiodic and is persistently recurrent, $P$ has at
least two successors and has only finitely many successors. Let
$\Gamma(P)$ be the {\it{last successor}} of $P$. Then there exists
an integer $q\geqslant 1$, largest among all of the successors of
$P$, such that $f^q(\Gamma(P))= P$.

We state some facts which will be used in the following as

\begin{enumerate}

\vspace{-0.2cm} \item[(F1)]\; $f^t (B(I))= I$, $\deg(f^t|_{B(I)})
\leqslant {d_{max}}^{b^2}$ and $\#\{i| \, c_0\in f^i (B(I)),
0\leqslant i <t\} \leqslant b$,

\vspace{-0.2cm}\item[(F2)]\; $f^s(A(I))= I$, $\deg(f^s|_{A(I)})
\leqslant {d_{max}}^{b^2+b}$ and $\#\{i| \, c_0\in f^i (A(I)),
0\leqslant i <s\} \leqslant b+1$,

\vspace{-0.2cm}\item[(F3)]\; $(B(I)-A(I)) \cap \textrm{orb}([c_0])
=\emptyset$,

\vspace{-0.2cm} \item[(F4)]\; $f^q(\Gamma(P))= P$ and
$\deg(f^q|_{\Gamma(P)})\leqslant {d_{max}}^{2b-1}$,

\vspace{-0.2cm} \item[(F5)]\; $f^i (\Gamma(P))$ does not contain
$c_0$ for all $0<i<q$.
\end{enumerate}

Now we can define the {\it{KSS nest}} in the following way: $I_0$
is a given piece containing $c_0$ and for $n\geqslant 0$,
\begin{eqnarray*}
&&L_n=A(I_n),\\
&&M_{n,0}=K_n=B(L_n),\\
&&M_{n,j+1}=\Gamma(M_{n,j}) \textrm{ for } 0\leqslant j\leqslant T-1,\\
&&I_{n+1}=M_{n,T}=\Gamma^T(K_n)=\Gamma^T(B(A(I_n))),
\end{eqnarray*}
with $T=3b$.

\begin{figure}[h]
\psfrag{In}{$I_n$}\psfrag{Tc0}{$T(c_0):$}

\psfrag{BIn}{$B(I_n)$} \psfrag{fsn}{$f^{s_n}$}

\psfrag{LnAIn}{$L_n=A(I_n)$} \psfrag{Knp}{$K_n'$}

\psfrag{ftn}{$f^{t_n}$} \psfrag{KnBLn}{$K_n=B(L_n)$}

\psfrag{gammaKn}{$\Gamma(K_n)$} \psfrag{fqn1}{$f^{q_{n,1}}$}

\psfrag{fqn}{$f^{q_n}$}
\psfrag{In1gammaTKn}{$I_{n+1}=\Gamma^T(K_n)$} \centering
\includegraphics[width=10cm]{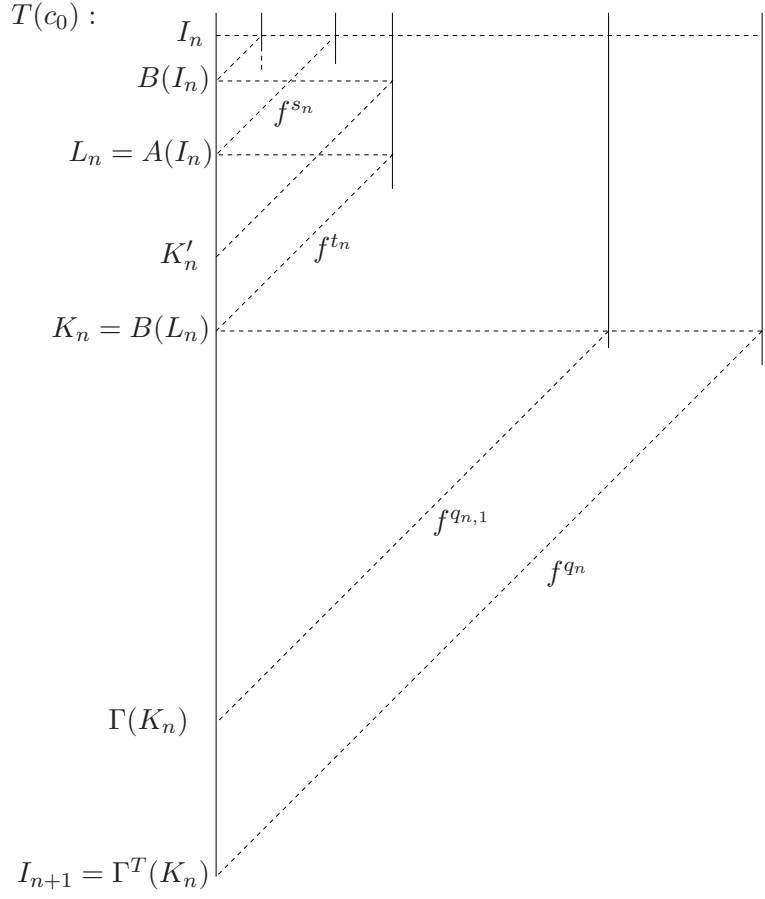} \caption{KSS
nest} \label{fig6}
\end{figure}

Suppose $f^{s_n}(L_n)= I_n$, $f^{t_n}(K_n)= L_n$,
$f^{q_{n,j}}(M_{n,j})= M_{n,j-1}$ for $1\leqslant j\leqslant T$,
and $q_n=\sum_{j=1}^T q_{n,j}.$ See Figure~\ref{fig6}.
\renewcommand{\captionlabeldelim}{~}

Let $p_n=q_{n-1}+s_n+t_n$. Then $f^{p_n}(K_n)= K_{n-1}$. From
(F1), (F2), and (F4), we have
$$d_0^{3b+2}\leqslant \deg(f^{p_n}|_{K_n}) \leqslant d_1,$$
where $d_1=d_{max}^{8b^2-2b}$.

For any puzzle piece $J$ containing $c_0$, let
$$r(J)=\min\{k(z)\,|\,\,z\in D(J)\cap J\},$$
where $k(z)$ is the smallest positive integer such that
$f^{k(z)}(z)\in J$. It is easy to prove that

\begin{enumerate}
\item[(1)] $r(J_1)\geqslant r(J_2)$ if $J_1\subset J_2$.
\item[(2)] $r(J)\geqslant k$ if $c_0\in J\subset J^\prime$, $f^k:\, J\to J^\prime$ and $c_0\not\in f^i (J)$ for $0<i<k$.
\end{enumerate}

The following lemma is a slight modification of Lemma 8.2 in
\cite{KSS}, and the proof is very much the same.

\begin{lemma}\label{lemma3.1}
Let $T=3b$. Then
\begin{enumerate}
\item[(1)] $r(I_n)\leqslant s_n \leqslant (b+1)r(L_n)$;
\item[(2)] $r(L_n)\leqslant t_n \leqslant br(K_n)$;
\item[(3)] $2r(M_{n,j-1})\leqslant q_{n,j}\leqslant r(M_{n,j})$ for $1\leqslant j\leqslant
T$;
\item[(4)] $p_{n+1}\geqslant 2p_n$, $\dfrac{p_n}{t_n}\leqslant b+4$.
\end{enumerate}
\end{lemma}

\begin{proof}
(1)\; The inequality $r(L_n)\leqslant s_n$ is obvious. Let
$G_j=f^j(L_n)$ and $0=j_0<j_1< \cdots <j_v=s_n$ be all the
integers such that $c_0\in G_{j_i}$. Then $v\leqslant b+1$ and
$f^{j_{i+1}-j_i}:G_{j_i}\to G_{j_{i+1}}$. Note that $c_0\notin
G_k$ for $j_i<k<j_{i+1}$. Hence $j_{i+1}-j_i\leqslant
r(G_{j_i})\leqslant r(G_0)=r(L_n)$ and
$$s_n=\sum_{i=0}^{v-1}(j_{i+1}-j_i)\leqslant vr(L_n)\leqslant (b+1)r(L_n). $$

(2)\; The proof is similar to that of (1).

(3)\; Let $1\leqslant j\leqslant T$. Since $M_{n,j}$ is a
successor of $M_{n,j-1}$ with $f^{q_{n,j}}:M_{n,j}\to M_{n,j-1}$
for some $q_{n,j}$ and $c_0\notin f^i(M_{n,j})$ for $0<i<q_{n,j}$,
we have $q_{n,j}\leqslant r(M_{n,j})$. Let $k$ be the smallest
positive integer with $f^k(c_0)\in M_{n,j-1}$ and
$J=L_{c_0}(M_{n,j-1})$. Then $f^k(J)=M_{n,j-1}$ and $J$ is the
first successor. Because $M_{n,j-1}$ has at least two successors
and $M_{n,j}$ is the last one, we have $q_{n,j}-k>0$. Denote
$x=f^k(c_0)$. Then $x\in M_{n,j-1}\cap D(M_{n,j-1})$ and
$f^{q_{n,j}-k}(x)\in M_{n,j-1}$. It follows that
$q_{n,j}=(q_{n,j}-k)+k\geqslant
r(M_{n,j-1})+r(M_{n,j-1})=2r(M_{n,j-1})$.

(4)\; By (3),
$$2^jr(K_n)=2^jr(M_{n,0})\leqslant q_{n,j}\leqslant
\frac{1}{2^{T-j}}r(M_{n,T})=\frac{1}{2^{T-j}}r(I_{n+1})$$ for any
$n\geqslant 1$ and $1\leqslant j\leqslant T$. From (1) and (2), we
have
\begin{align*}
p_{n+1}&=q_n+s_{n+1}+t_{n+1}\\
&=\sum_{j=1}^Tq_{n,j}+s_{n+1}+t_{n+1}\\
&\geqslant (2^{T+1}-2)r(K_n)+r(I_{n+1})+r(L_{n+1})\\
&\geqslant 2^{T+1}r(K_n)=2^{3b+1}r(K_n)
\end{align*}
and
\begin{align*}
p_n&=q_{n-1}+s_n+t_n\\
&=\sum_{j=1}^Tq_{n-1,j}+s_n+t_n\\
&< 2r(I_n)+(b+1)r(L_n)+br(K_n)\\
&\leqslant (2b+3)r(K_n).
\end{align*}
Therefore, $p_{n+1}\geqslant 2p_n$.

The second inequality can be obtained from the following fact
\begin{align*}
p_n&=q_{n-1}+s_n+t_n\\
&< 2r(I_n)+(b+1)r(L_n)+t_n\\
&\leqslant (b+3)r(L_n)+t_n\\
&\leqslant (b+4)t_n.
\end{align*}
\end{proof}

\section{Proof of the Main Proposition}

\noindent Let $K_n^\prime= \textrm{Comp}_{c_0} f^{-t_n} (B(I_n))$.
The conditions $F(c_0)=[c_0]$ and $(B(I_n)-L_n) \cap
\textrm{orb}([c_0])=\emptyset$ imply
$$d_0\leqslant \deg(f^{t_n}|_{K_n^\prime}) =\deg(f^{t_n}|_{K_n}) \leqslant d_1. $$
Let $\mu_n=\!\!\!\mod(K_n^\prime - \overline{K_n})$.

The main result in this section is the following lemma.
\begin{lemma}\label{lemma4.1}
$\liminf_{n\to\infty} \mu_n >0$.
\end{lemma}

We first state a covering lemma recently given by Kahn and Lyubich
which will play a crucial rule in the proof of Lemma
\ref{lemma4.1}.

\begin{KLlemma}[\cite{KL-1}]
Fix some $\eta>0$. Let $A\subset\subset A^\prime \subset$
$\mathrm{int}$$(U)$ and $B\subset\subset B^\prime\subset$
$\mathrm{int}$$(V)$ be two nests of Jordan disks. Let $f:\,
(U,A^\prime,A)\to (V,B^\prime, B)$ be a holomorphic proper mapping
between the respective disks, and let $D=\deg(f|_U)$ and
$d=\deg(f|_{A^\prime})$. Assume the following collar property
$$\!\!\!\mod(B^\prime-\overline{B})\geqslant \eta \!\!\!\mod (U-\overline{A}).$$
Then there exists an $\epsilon>0$ (depending on $\eta$ and $D$) such
that
$$\!\!\!\mod(V-\overline{B})\leqslant C\eta^{-1}d^2 \!\!\!\mod (U-\overline{A})$$
or
$$\!\!\!\mod(U-\overline{A})\geqslant \epsilon,$$
where $C$ is an absolute constant.
\end{KLlemma}

\noindent{\it Proof of Lemma 3}

Suppose $\liminf_{n\to\infty} \mu_n=0$. Let
$\mu_{k_n}=\min\{\mu_1,\mu_2,\ldots,\mu_n\}$. Then
$\lim_{n\to\infty}k_n=\infty$ and $\lim_{n\to\infty}\mu_{k_n}=0$.
Take an integer $j_0$ satisfying
$$ 2^{3b(j_0-1)} \geqslant (b+1)(2b+9)$$
and a large integer $N$. Let
$M=p_{k_n-j_0}+p_{k_n-j_0-1}+\cdots+p_{k_n-N+1}$. Then
\begin{enumerate}
\item[(1)] $M<2p_{k_n-j_0}$ (by Lemma 2(4)),
\item[(2)] $f^M(K_{k_n-j_0})= K_{k_n-N}$,
\item[(3)] $d_0^{(3b+2)(N-j_0)}\leqslant D = \deg(f^M|_{K_{k_n-j_0}})\leqslant {d_1}^{N-j_0}$,
\end{enumerate}
where $d_0=\deg_{c_0}f$ and $d_1$ is the constant in section 3
depending only on $b$ and $d_{max}$.

\begin{figure}[h]
\psfrag{fM}{$f^M$}

\psfrag{x}{\tiny$x$} \psfrag{Ax}{\tiny$A_x$}
\psfrag{Axp}{\tiny$A_x'$} \psfrag{Kknj0}{\small $K_{k_n-j_0}$}

\psfrag{fl}{\tiny$f^l$} \psfrag{Kknj0p}{\small $K_{k_n-j_0}'$}

\psfrag{y}{\tiny$y$} \psfrag{By}{\tiny$B_y$}

\psfrag{Byp}{\tiny$B_y'$} \psfrag{KknN}{$K_{k_n-N}$}

\centering
\includegraphics[width=11cm]{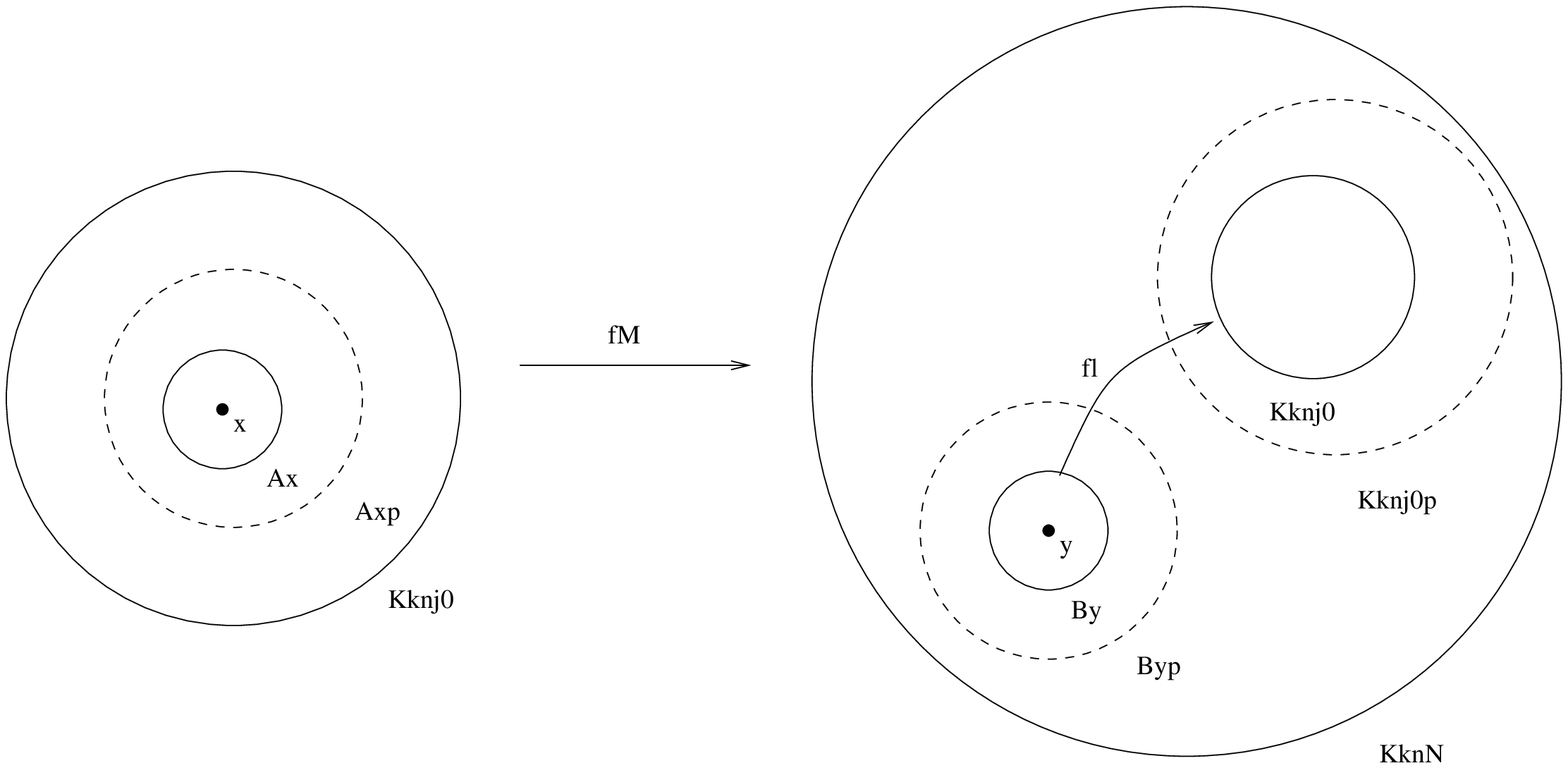}\caption{}\label{fig7}
\end{figure}

For any $x\in K_{k_n-j_0}\cap \textrm{orb}([c_0])$, let $y=f^M
(x)$, $B_y=\hat{L}_y(K_{k_n-j_0})$, $f^l(B_y)=K_{k_n-j_0}$ and
$B_y^\prime=\hat{L}_y(K_{k_n-j_0}^\prime)$. See Figure~\ref{fig7}.

Let $A_x=\textrm{Comp}_x f^{-M}(B_y)$ and
$A_x^\prime=\textrm{Comp}_x f^{-M}(B_y^\prime)$. From the conditions
$(K_{k_n-j_0}^\prime - K_{k_n-j_0}) \cap
\textrm{orb}([c_0])=\emptyset$ and (T3), we have
$$ \deg(f^M|_{A_x^\prime})  =  \deg(f^M |_{A_x}).$$

\begin{figure}[h]
\psfrag{Tx}{$T(x):$} \psfrag{Kknj0Pn0c0}{$K_{k_n-j_0}$}

\psfrag{ftknj0}{\tiny$f^{t_{k_n-j_0}}$} \psfrag{W1}{$W_1$}
\psfrag{n00}{\tiny
$(n_0,0)$}\psfrag{n0m1}{\tiny$(n_0,m_1)$}\psfrag{n0m2}{\tiny$(n_0,m_2)$}\psfrag{n0m}{\tiny$(n_0,M)$}
\psfrag{fl1}{\tiny$f^{l_1}$} \psfrag{W2}{$W_2$}
\psfrag{By}{$B_y$}\psfrag{fM}{$f^M$}\psfrag{fl}{$f^l$}
\psfrag{n0mL}{\tiny$(n_0,m_L)$} \psfrag{fl2}{\tiny$f^{l_2}$}
\psfrag{n0Ml}{\tiny$(n_0,M+l)$}\psfrag{Ax}{$A_x$} \centering
\includegraphics[width=10cm]{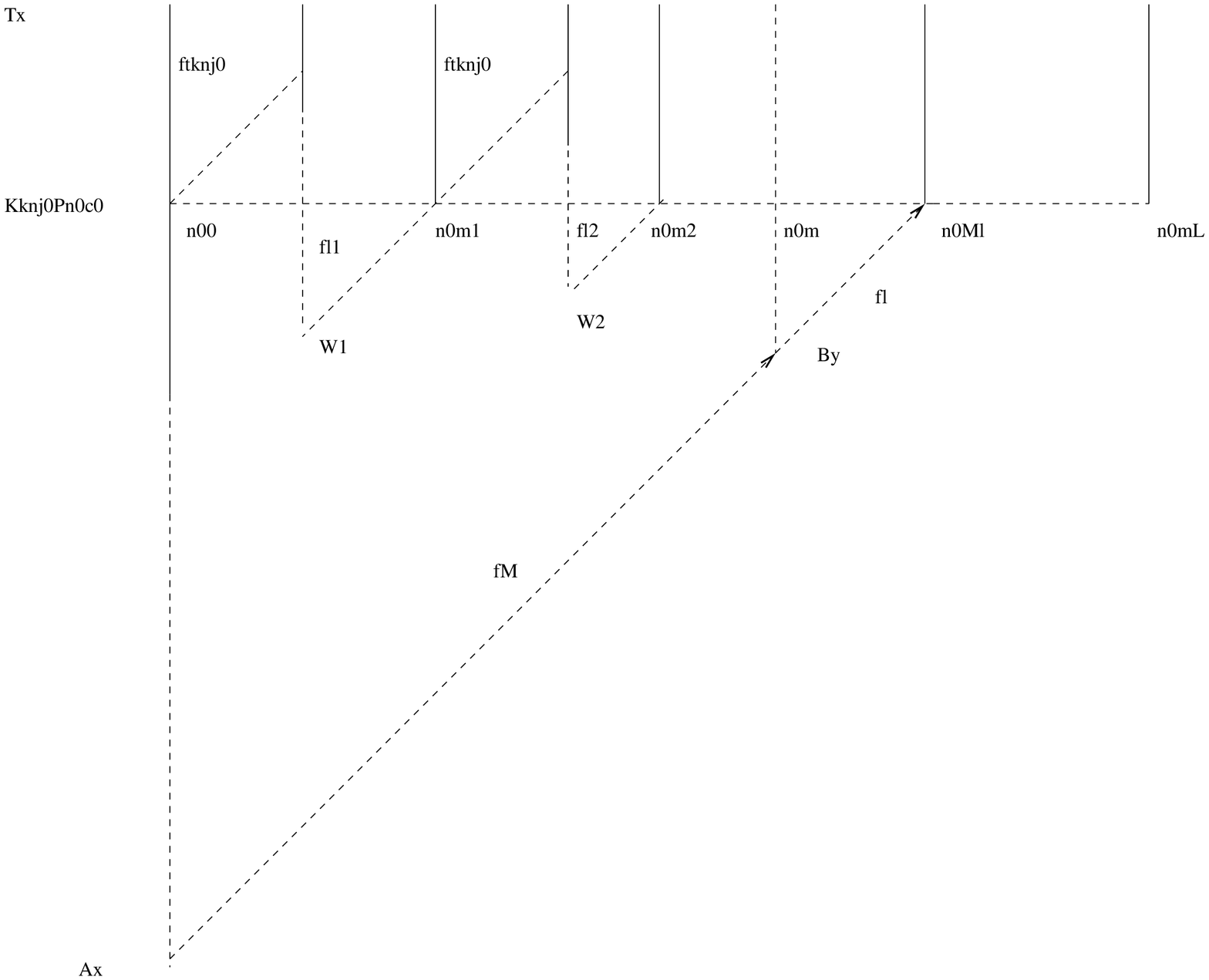}\caption{}\label{fig8}
\end{figure}

\begin{claim1}
For any $x\in K_{k_n-j_0}\cap \mathrm{orb}([c_0])$,
$$\deg(f^M|_{A_x^\prime})  =  \deg(f^M |_{A_x})\leqslant d_2,$$
where $d_2$ is a constant depending only on $b$ and $d_{max}$.
\end{claim1}

\begin{proof}
Let $K_{k_n-j_0}=P_{n_0}(c_0)$. Suppose $(n_0,m_1)$ is the first
$c_0$-position on the right of $(n_0,t_{k_n-j_0})$ in $T(x)$. Let
$W_1= P_{n_0+m_1-t_{k_n-j_0},t_{k_n-j_0}}(x)$ and $l_1=
m_1-t_{k_n-j_0}$. Then $f^{l_1}(W_1)= P_{n_0}(c_0)$ and $\deg
(f^{l_1}|_{W_1})\leqslant d_{max}^b$. Repeating this process, we
have infinitely many $c_0$-positions $\{(n_0,m_i)\}_{i\geqslant
1}$ such that $(n_0,m_i)$ is the first $c_0$-position on the right
of $(n_0,m_{i-1}+t_{k_n-j_0})$ for each $i\geqslant 1$ in $T(x)$.
For any $i\geqslant 1$, $f^{l_i}(W_i)= P_{n_0}(c_0)$ and $\deg
(f^{l_i}|_{W_i})\leqslant d_{max}^b$, where
$W_i=P_{n_0+l_i,m_{i-1}+t_{k_n-j_0}}(x)$,
$l_i=m_i-m_{i-1}-t_{k_n-j_0}$. Let $L\geqslant 1$ be the smallest
integer such that $m_L\geqslant M+l$. Then
$(L-1)t_{k_n-j_0}\leqslant m_{L-1}<M$. By Lemma 2(4) and (F1),
$$L-1\leqslant \dfrac{M}{t_{k_n-j_0}} <
\dfrac{2p_{k_n-j_0}}{t_{k_n-j_0}}\leqslant 2(b+4)=2b+8$$ and
$$\deg(f^M|_{A_x}) \leqslant
\left(\deg(f^{t_{k_n-j_0}}|_{K_{k_n-j_0}})\right)^L\cdot
\prod_{j=1}^L (\deg f^{l_i}|_{W_i}) \leqslant d_2,$$ where
$d_2=d_{max}^{(b^2+b)(2b+9)}$. See Figure~\ref{fig8}.
\end{proof}

Suppose $f^{M^\prime}(I_{k_n})= K_{k_n-j_0}$ and
 $f^\sigma(B(I_{k_n}))= I_{k_n}$, where $$M^\prime =
q_{k_n-1}+p_{k_n-1}+\cdots + p_{k_n-j_0+1}.$$ From (F1) and
$\deg(f^{p_n}|_{K_n})\leqslant d_1$, we have
$\deg(f^\sigma|_{B(I_{k_n}})\leqslant d_{max}^{b^2}$ and
$$\deg(f^{M^\prime}|_{I_{k_n}})\leqslant d_{max}^{b^2}d_1^{j_0-1},$$
where $d_1=d_{max}^{8b^2-2b}$ is obtained in section 3.

\begin{figure}[h]
\psfrag{Lkn}{$L_{k_n}$} \psfrag{BIkn}{$B(I_{k_n})$}
\psfrag{c0}{\tiny$c_0$} \psfrag{fsigma}{$f^{\sigma}$}
\psfrag{fMp}{$f^{M'}$}\psfrag{x}{$x$}\psfrag{Ax}{$A_x$}

\psfrag{Omega}{$\Omega$}
\psfrag{Ikn}{$I_{k_n}$}\psfrag{Kknj0}{$K_{k_n-j_0}$}

\centering\includegraphics[width=9cm]{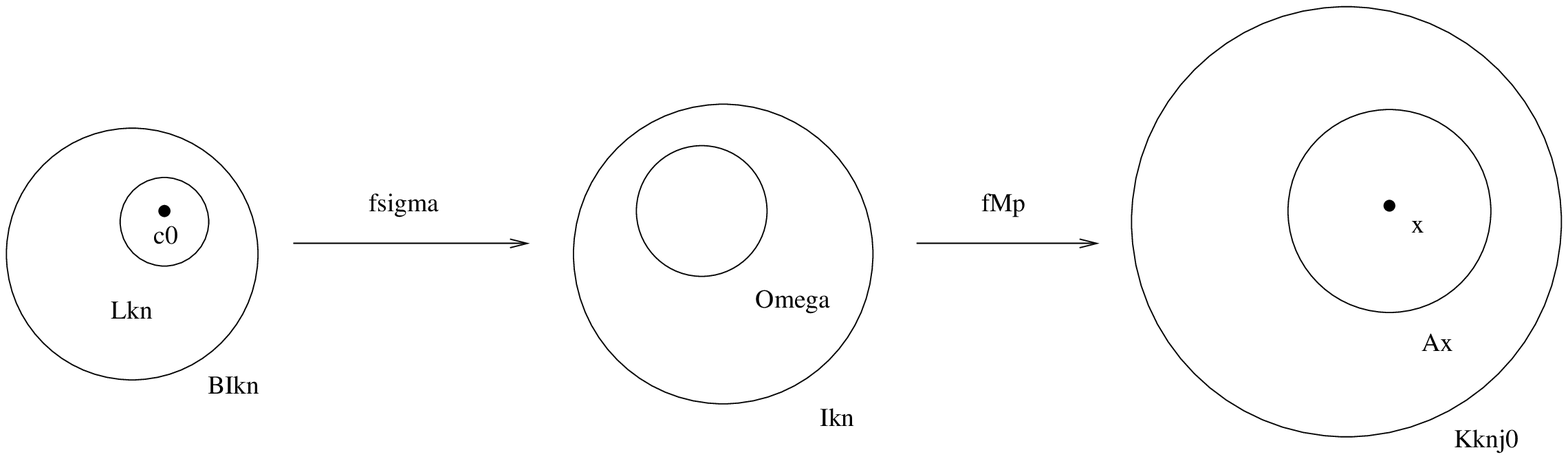}\caption{}\label{fig9}
\end{figure}

Let $x= f^{M^\prime+\sigma} (c_0)$ and let $A_x$ be the puzzle
piece constructed as above. See Figure \ref{fig9}.

\begin{figure}[h]
\psfrag{fMl}{$f^{M+l}$}\psfrag{Omega}{$\Omega$}
 \psfrag{x}{$x$}\psfrag{Ax}{$A_x$}\psfrag{Kknj0}{$K_{k_n-j_0}$}
\psfrag{fMp}{$f^{M'}$} \psfrag{fMpOmega}{$f^{M'}(\Omega)$}
\psfrag{fr}{$f^r$}
\centering\includegraphics[width=9cm]{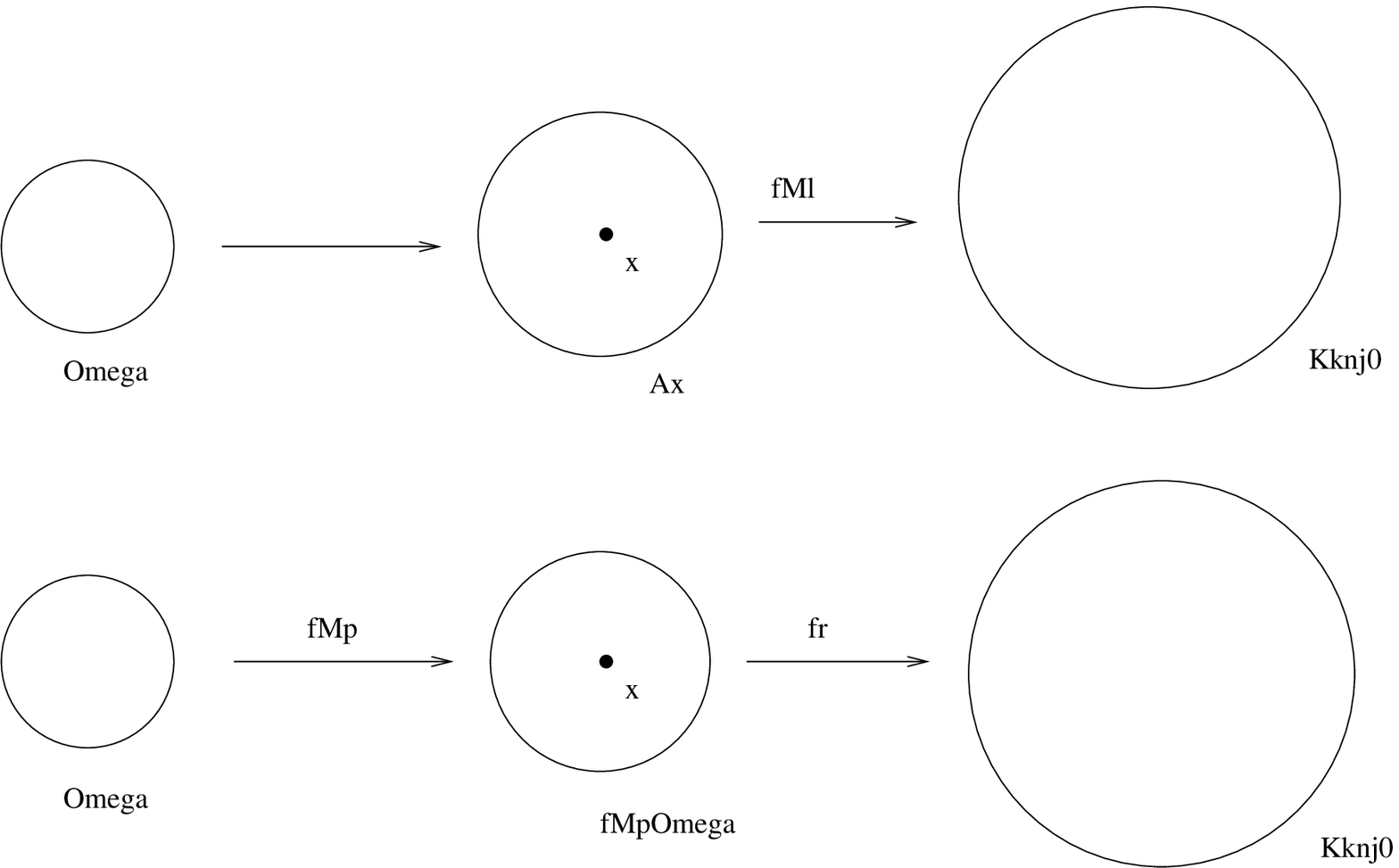}\caption{}\label{fig10}

\end{figure}
\begin{figure}[h]
\psfrag{n0plus}{\tiny$(n_0+v_{i+1}-v_i,v_i)$}\psfrag{n0plusk}{\tiny$(n_0+k_i,v_{i+1}-k_i)$}\psfrag{c}{$c$}
\psfrag{fr}{$f^{r}$}\psfrag{n0mi}{$(n_0,v_i)$}\psfrag{n0mpsigma}{\small$(n_0,M'+\sigma)$}\psfrag{n0mpsigmar}
{$(n_0,M'+\sigma+r)$}
 \psfrag{n0mi1}{$(n_0,v_{i+1})$}\psfrag{Tx}{$T(c_0):$}\psfrag{fmpsigma}{$f^{M'+\sigma}$}
\psfrag{c0}{$c_0$} \psfrag{fMpOmega}{$f^{M'}(\Omega)$}
\psfrag{fr}{$f^r$}\psfrag{Lkn}{$L_{k_n}$} \psfrag{Kknj0Pn0c0}{
$\begin{array}{lll} P_{n_0}(c_0)&\ &\\
 =K_{k_n-j_0}&\ &
\end{array}
$}
\centering\includegraphics[width=14cm]{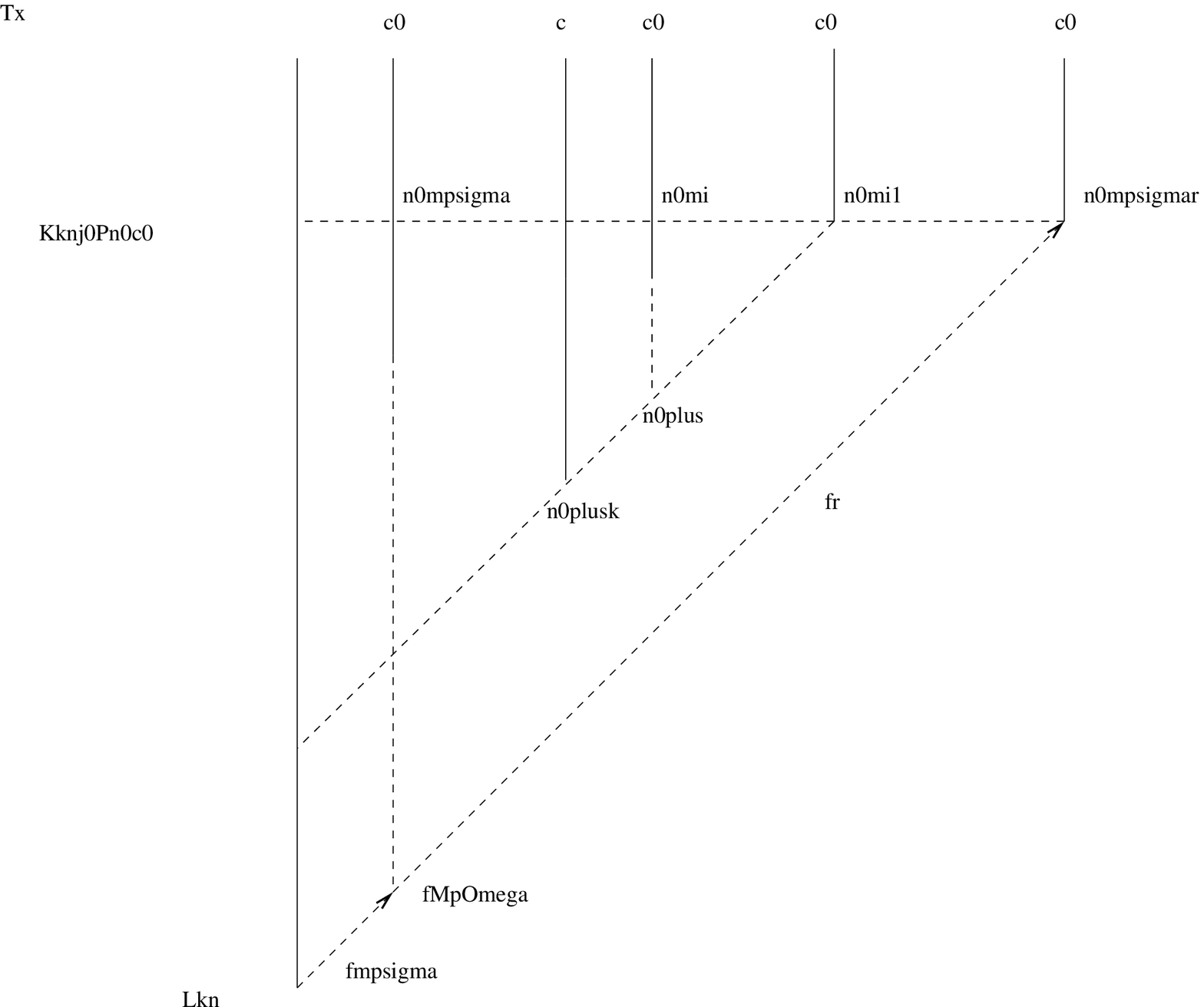}\caption{}\label{fig11}
\end{figure}

\begin{claim2}
Let $\Omega=f^\sigma (L_{k_n})\subset I_{k_n}$. Then
$f^{M^\prime}(\Omega)\subset A_x$.
\end{claim2}

\begin{proof}
Suppose $f^r(\Omega)= I_{k_n}$, then $f^{r+M^\prime}(\Omega)=
K_{k_n-j_0}$ and $r\geqslant r(I_{k_n})$. See Figure \ref{fig10}.

Let $K_{k_n-j_0}=P_{n_0}(c_0)$. Suppose
$(n_0,v_0),(n_0,v_1),\ldots,(n_0,v_k)$ are all $c_0$-positions
between $(n_0,M^\prime+\sigma)$ and $(n_0,M^\prime+\sigma+r)$ in
$T(c_0)$ with $v_0=M^\prime + \sigma$ and $v_k=M^\prime + \sigma
+r$. See Figure~\ref{fig11}.

\begin{subclaim}
For all $0\leqslant i\leqslant k-1$, $v_{i+1}-v_i\leqslant
q_{k_n-j_0,1}$.
\end{subclaim}
\begin{proof}
We recall that $q_{k_n-j_0,1}\geqslant 1$ is the integer, largest
among all of the successors of $P_{n_0}(c_0)$, such that
$f^{q_{k_n-j_0,1}}(\Gamma(P_{n_0}(c_0)))= P_{n_0}(c_0)$.

If $(n_0+v_{i+1}-v_i,v_i)$ is critical, then
$P_{n_0+v_{i+1}-v_i}(c_0)$ is a successor of $P_{n_0}(c_0)$ and
$v_{i+1}-v_i\leqslant q_{k_n-j_0,1}$.

Suppose that $(n_0+v_{i+1}-v_i,v_i)$ is not critical. Let $k_i$ be
the smallest integer between $v_{i+1}-v_i$ and $v_{i+1}$ such that
$(n_0+k_i,v_{i+1}-k_i)$ is a critical position, say it is a
$c$-position, see Figure~\ref{fig11}. Then $c\in [c_0].$ If
$c=c_0$, then $P_{n_0+k_i}(c_0)$ is a successor of $P_{n_0}(c_0)$
and $v_{i+1}-v_i < k_i\leqslant q_{k_n-j_0,1}$. If $c\neq c_0$,
let $P_{n_0+l_i}(c_0)=L_{c_0}(P_{n_0+k_i}(c))$. Then
$P_{n_0+l_i}(c_0)$ is a successor of $P_{n_0}(c_0)$ and
$v_{i+1}-v_i < k_i<l_i\leqslant q_{k_n-j_0,1}$.

\end{proof}
By the Subclaim and Lemma 2(3),
$$\#\{i|\,f^i(f^{M^\prime}(\Omega))\subset K_{k_n-j_0},0\leqslant i<r\}
\geqslant \dfrac{r}{q_{k_n-j_0,1}}\geqslant
\dfrac{r(I_{k_n})}{r(I_{k_n-j_0+1})}\geqslant 2^{3b(j_0-1)},$$ since
$r(I_{n+1})\geqslant 2^{3b}r(K_n)\geqslant 2^{3b}r(I_n)$ for all
$n\geqslant 0$. See Figure~\ref{fig11}.

By $K_{k_n-j_0}=B(L_{k_n-j_0})$ and (F1), we have
$$\#\{i|\,c_0\in
f^i(K_{k_n-j_0}), 0\leqslant i\leqslant t_{k_n-j_0}\}\leqslant
b+1.$$ For each $1\leqslant j \leqslant L$, $(n_0,k)$ is not
$c_0$-position for $m_j-l_j< k < m_j$. Hence
\begin{align*}
\#\{i|\,f^i(A_x)\subset K_{k_n-j_0}, 0\leqslant i<M+l\} &\leqslant
L\cdot \#\{i|\,c_0\in f^i(K_{k_n-j_0}), 0\leqslant i\leqslant
t_{k_n-j_0}\}\\
&\leqslant (b+1)L <(b+1)(2b+9).
\end{align*}
The integers $L$, $m_j$, and $l_j$ are the same as in the proof of
Claim 1. See Figure~\ref{fig8}.

The condition $2^{3b(j_0-1)}\geqslant (b+1)(2b+9)$ implies that
$$f^{M^\prime}(\Omega) \subset A_x.$$
\end{proof}

\vskip0.5cm

\begin{claim3}
There exists a positive constant $\eta$ depending only on $b$ and
$d_{max}$ such that
$$\!\!\!\mod(B_y^\prime-\overline{B_y}) \geqslant \eta \!\!\!\mod(K_{k_n-j_0}-\overline{A_x}).$$
\end{claim3}

\begin{proof}
Since $\deg(f^l|_{B_y}) = \deg (f^l|_{B_y^\prime})\leqslant
d_{max}^{b-1}$, we have
\begin{eqnarray*}
\!\!\!\mod(B_y^\prime-\overline{B_y})\hspace{-0.6cm}&&=
\dfrac{1}{\deg(f^l|_{B_y})}\!\!\!\mod(K_{k_n-j_0}^\prime - \overline{K_{k_n-j_0}}) \\
&&\geqslant {d_{max}}^{-(b-1)} \mu_{k_n-j_0} \\
&&\geqslant {d_{max}}^{-(b-1)} \mu_{k_n}.
\end{eqnarray*}

\begin{figure}[h]
\psfrag{Kkn}{\tiny$K_{k_n}$} \psfrag{Kknp}{\tiny$K'_{k_n}$}
\psfrag{c0}{\tiny$c_0$}
\psfrag{ftkn}{\tiny$f^{t_{k_n}}$}\psfrag{Lkn}{\tiny$L_{k_n}$}\psfrag{fsigma}{\tiny$f^{\sigma}$}\psfrag{Omega}{\tiny$\Omega$}
\psfrag{Ikn}{\tiny$I_{k_n}$}\psfrag{fMpOmega}{\tiny$f^{M'}(\Omega)$}\psfrag{Ax}{\tiny$A_x$}\psfrag{fMp}{\tiny$f^{M'}$}
\psfrag{Kknj0}{\tiny$K_{k_n-j_0}$}\psfrag{BIkn}{\tiny$B(I_{k_n})$}
\psfrag{x}{\tiny$x$}
\centering\includegraphics[width=13cm]{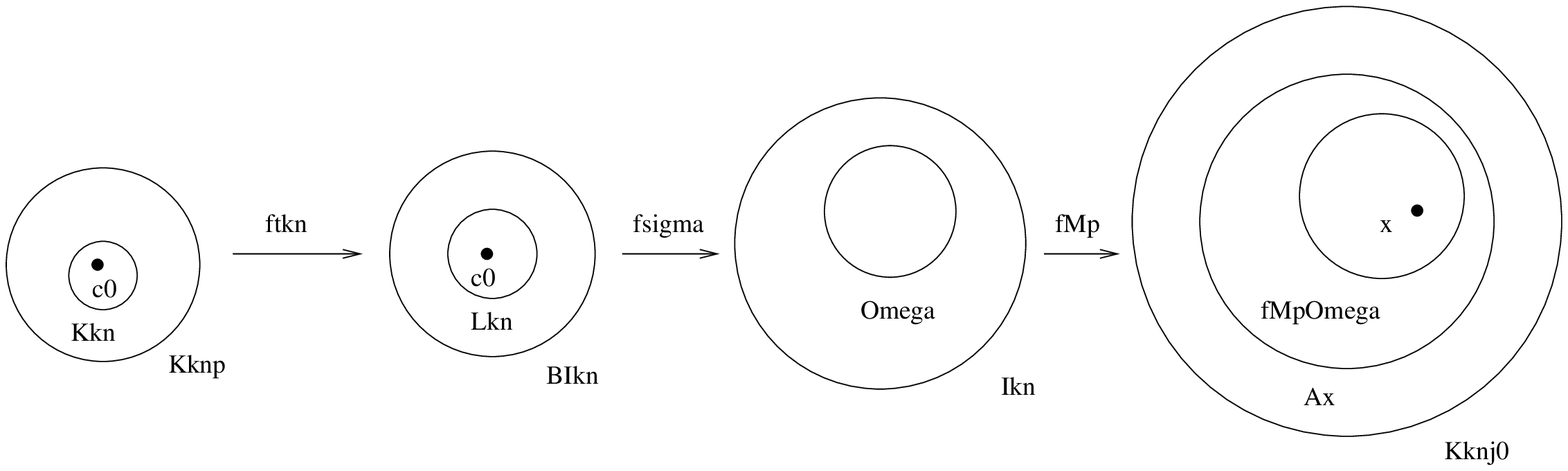}\caption{}\label{fig12}
\end{figure}

\vskip0.5cm By Claim 2 and (F1) in section 3,
\begin{eqnarray*}
\mod(K_{k_n-j_0}-\overline{A_x})\hspace{-0.6cm}&&\leqslant
\!\mod(K_{k_n-j_0}-\overline{f^{M^\prime} (\Omega)}) \\
&&\leqslant \deg(f^{t_{k_n}+\sigma+M^\prime}|_{K_{k_n}^\prime})\cdot\mod(K_{k_n}^\prime - \overline{K_{k_n}}) \\
&&=\deg(f^{t_{k_n}}|_{K_{k_n}^\prime})\deg(f^\sigma|_{B(I_{k_n}})\deg(f^{M^\prime}|_{I_{k_n}})\cdot\mu_{k_n}\\
&&=\deg(f^{t_{k_n}}|_{K_{k_n}})\deg(f^\sigma|_{B(I_{k_n}})\deg(f^{M^\prime}|_{I_{k_n}})\cdot\mu_{k_n}\\
&&\leqslant d_{max}^{b^2}\cdot d_{max}^{b^2}\cdot d_{max}^{b^2}d_1^{j_0-1} \cdot\mu_{k_n}\\
&&= d_3\mu_{k_n},
\end{eqnarray*}
where $d_3=d_{max}^{3b^2}d_1^{j_0-1}$ is a constant depending only
on $b$ and $d_{max}$. See Figure~\ref{fig12}.

Take $\eta=d_3^{-1}{d_{max}}^{-(b-1)}$. We have
$$\!\!\!\mod(B_y^\prime - \overline{B_y}) \geqslant \eta
\!\!\!\mod(K_{k_n-j_0}-\overline{A_x}).$$
\end{proof}

\begin{figure}[h]
\psfrag{KknN}{\tiny$K_{k_n-N}$} \psfrag{x}{\tiny$x$}
\psfrag{Ax}{\tiny$A_x$} \psfrag{Axp}{\tiny$A_x'$}
\psfrag{Kknj0}{\tiny$K_{k_n-j_0}$}\psfrag{y}{\tiny$y$}
\psfrag{By}{\tiny$B_y$} \psfrag{fM}{\tiny$f^M$}
\psfrag{Byp}{\tiny$B_y'$}
\centering\includegraphics[width=10cm]{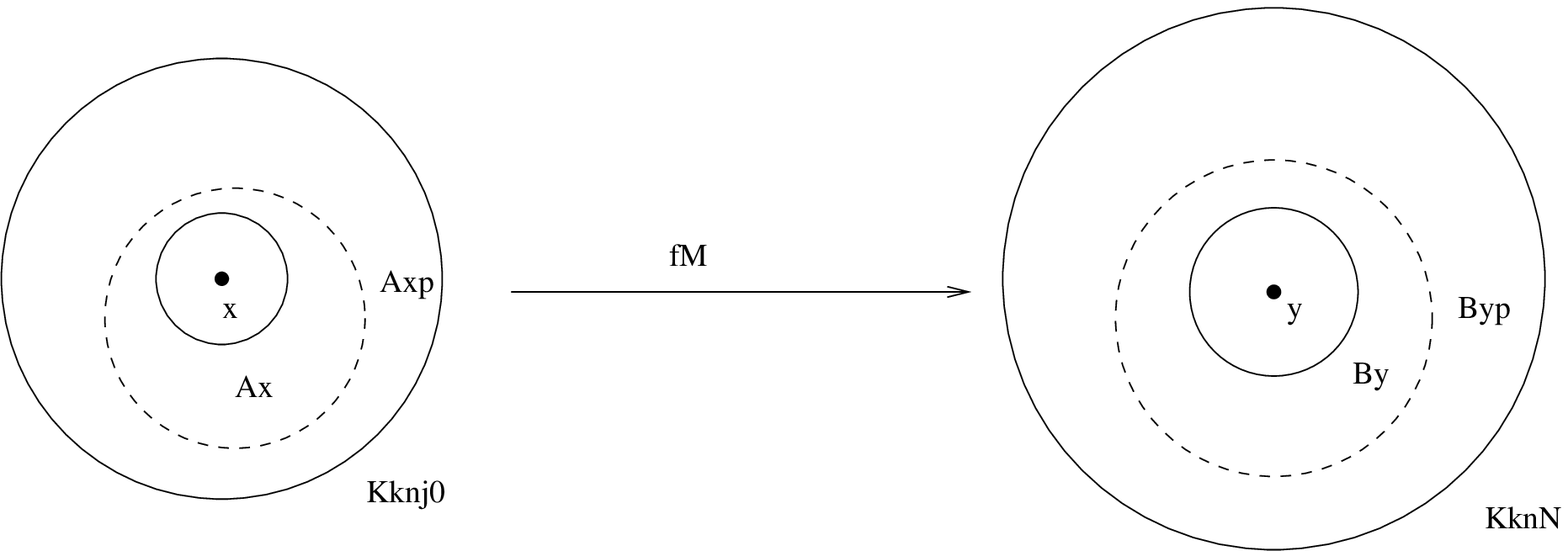}\caption{}\label{fig13}
\end{figure}

Now we have a holomorphic proper mapping $f^{M}:\, (K_{k_n-j_0},
A_x^\prime, A_x) \to (K_{k_n-N},B_y^\prime,B_y)$ satisfying
\begin{enumerate}
\item[(1)] $d_0^{(3b+2)(N-j_0)} \leqslant D = \deg( f^M|_{K_{k_n-j_0}})
\leqslant d_1^{N-j_0}$,
\item[(2)] $\deg(f^M|_{A_x^\prime}) = \deg(f^M|_{A_x}) \leqslant d_2$,
\item[(3)] $\!\!\!\mod(B_y^\prime-\overline{B_y}) \geqslant \eta \!\!\!\mod(K_{k_n-j_0}-\overline{A_x})$,
\end{enumerate}
where $d_0=\deg_{c_0}f$ and $d_1,d_2,\eta$ are constants depending
only on $b$ and $d_{max}$. See Figure~\ref{fig13}.

By the Kahn-Lyubich Lemma,
\begin{equation}\label{KLles}
\!\!\!\mod(K_{k_n-N}-\overline{B_y}) \leqslant C\eta^{-1} d_2^2
\!\!\!\mod (K_{k_n-j_0}-\overline{A_x})
\end{equation}
or
\begin{equation}\label{KLgeq}
\!\!\!\mod(K_{k_n-j_0}-\overline{A_x})\geqslant \epsilon.
\end{equation}

\begin{figure}[h]
\psfrag{KknN}{$K_{k_n-N}$} \psfrag{Ty}{$T(y):$}
\psfrag{Kkni}{$K_{k_n-i}$}\psfrag{Kknj0}{$K_{k_n-j_0}$}
\psfrag{By}{$B_y$} \psfrag{Kknip}{$K_{k_n-i}'$}
\psfrag{fri}{$f^{r_i}$}

\psfrag{Viyp}{$V'_i(y)$} \psfrag{Viy}{$V_i(y)$} \centering
\includegraphics[width=10cm]{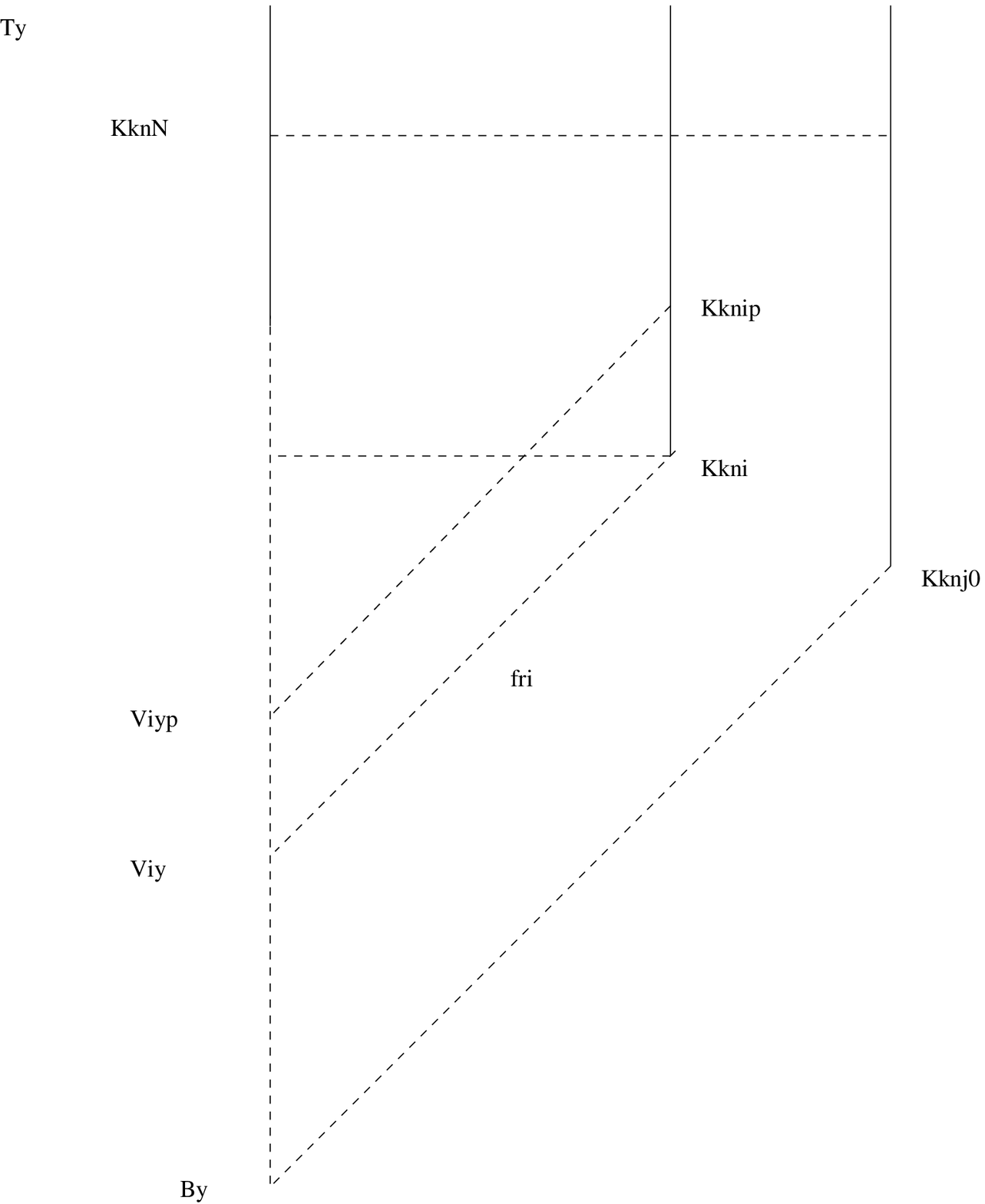}\caption{}\label{fig14}
\end{figure}

We first prove that the inequality (\ref{KLles}) is impossible for
$N$ large enough. For each $j_0\leqslant i\leqslant N-1$, let
$V_i(y)=L_y(K_{k_n-i})$, $V_i^\prime(y)=L_y(K_{k_n-i}^\prime)$,
and $f^{r_i}(V_i(y))= K_{k_n-i}$. See Figure~\ref{fig14}. Then
$f^{r_i}(V_i^\prime(y))= K_{k_n-i}^\prime$ and
$$\deg(f^{r_i}|_{V_i(y)}) = \deg(f^{r_i}|_{V_i^\prime(y)}) \leqslant {d_{max}}^b.$$

Therefore
\begin{eqnarray*}
\!\!\!\mod(V_i^\prime(y) - \overline{V_i(y)}) \hspace{-0.6cm}&&=
\dfrac{1}{\deg(f^{r_i}|_{V_i(y)})}
\!\!\!\mod(K_{k_n-i}^\prime - \overline{K_{k_n-i}}) \\
&&\geqslant {d_{max}}^{-b} \mu_{k_n-i} \\
&&\geqslant {d_{max}}^{-b} \mu_{k_n}
\end{eqnarray*}
and
$$\!\!\!\mod(K_{k_n-N}-\overline{B_y}) \geqslant (N-j_0) {d_{max}}^{-b} \mu_{k_n}. $$
By the proof of Claim 3,
$$ \!\!\!\mod(K_{k_n-j_0}-\overline{A_x}) \leqslant d_3 \mu_{k_n}.$$
Hence
$$ \!\!\!\mod(K_{k_n-N}-\overline{B_y})\geqslant (N-j_0)d_3^{-1} d_{max}^{-b} \!\!\!\mod (K_{k_n-j_0}-\overline{A_x}).$$
This implies that the inequality (\ref{KLles}) is impossible for
$N$ large enough.

Take a large $N_0$ such that (\ref{KLles}) does not hold. We have
$$\!\!\!\mod(K_{k_n-j_0}-\overline{A_x})\geqslant \epsilon>0,$$
where $\epsilon$ depends only on $\eta$ and $N_0$. This contradicts
the fact
$$\!\!\!\mod(K_{k_n-j_0}-\overline{A_x}) \leqslant d_3\mu_{k_n} \to 0$$
as $n\to \infty$. This completes the proof of Lemma 3. \qed

\vskip2.5cm

\noindent{\it Proof of the Main Proposition}

By Lemma 3, $\mu_{k_n}\geqslant \mu >0$ for some constant $\mu$. The
Gr\"{o}tzsch's inequality implies that
$$ \!\!\!\mod(P_0(c_0)-K_f(c_0)) \geqslant \sum_{k=0}^\infty \!\!\!\mod (K_{k_n}^\prime - \overline{K_{k_n}})=+\infty $$
and $K_f(c_0) = \bigcap_{n\geqslant 0} P_n(c_0) = \{c_0\}$.

For any $x\in K_f$ with $x\to c_0$, let $V_n(x)=L_x(K_{k_n})$,
$f^{r_n}(V_n(x))= K_{k_n}$, and let $V_n^\prime(x) =
\textrm{Comp}_x f^{-r_n}(K_{k_n}^\prime)$. The degree of
$f^{r_n}:\, V_n(x) \to K_{k_n}$ is uniformly bounded for all $n$.

If there are infinitely many $n$, say $\{n_j\}$, such that there
is at most one piece in
$\{V_n^\prime(x),f(V_n^\prime(x)),\ldots,f^{r_n}(V_n^\prime(x))=K_{k_n}^\prime\}$
containing $c$ for any $c\in \mathrm{Crit}-[c_0]$, then the degree
of $f^{r_{n_j}}:\, V_{n_j}^\prime(x) \to K_{k_{n_j}}^\prime$ is
uniformly bounded for all $j$ and there is a constant
$\widetilde{\mu}>0$ such that
$$\!\!\!\mod(V_{n_j}^\prime(x) - \overline{V_{n_j}(x)}) \geqslant \widetilde{\mu}$$
for all $j$. In this case, $K_f(x)=\{x\}$.

\begin{figure}[h]

\psfrag{Tx}{$T(x):$} \psfrag{c1}{$c_1$} \psfrag{c0}{$c_0$}
\psfrag{kknj}{$K_{k_{n_j}}$}

\psfrag{Kknjp}{$K_{k_{n_j}}'$} \psfrag{vnjx}{$V_{n_j}(x)$}

\psfrag{vnjxp}{$V'_{n_j}(x)$} \centering
\includegraphics[width=10cm]{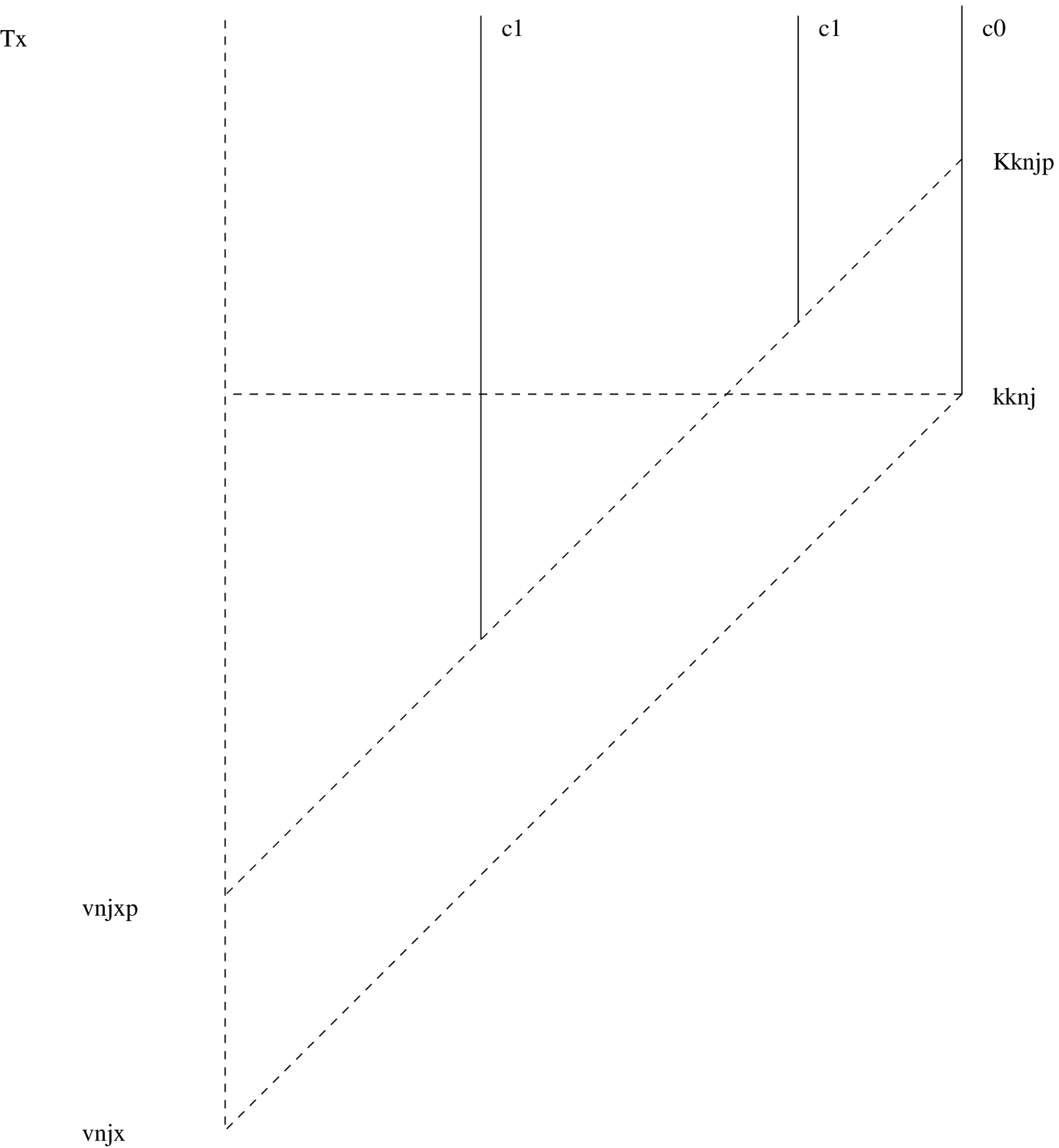}\caption{}\label{fig15}
\end{figure}

Suppose for each large $n$, there are two pieces in
$$\{V_n^\prime(x),f(V_n^\prime(x)), \ldots, f^{r_n}(V_n^\prime (x)) =
K_{k_n}^{\prime}\}$$ containing $c$ for some $c\in \mathrm{Crit} -
[c_0]$. There exist $c_1\in \mathrm{Crit}- [c_0]$ and a
subsequence $\{n_j\}$ such that there are two pieces in
$$\{V^\prime_{n_j}(x),
f(V_{n_j}^\prime(x)), \ldots,
f^{r_{n_j}}(V_{n_j}^\prime(x))=K_{k_{n_j}}^\prime\}$$ containing
$c_1$. See Figure~\ref{fig15}.

We conclude that $c_1\to c_1$ and $c_1\to c_0$. From the conditions
$(K_{k_{n_j}}^\prime - K_{k_{n_j}}) \cap
\textrm{orb}([c_0])=\emptyset$ and (T3), we have $c_1\not\in [c_0]$.
The Lemma 1 implies $T(c_1)$ is reluctantly recurrent. By
Proposition 1, $K_f(x)=\{x\}$ since $x\to c_1$.

This ends the proof of the Main Proposition. \qed

\section{Wandering components of filled-in Julia sets}

\noindent In this section, let $f$ be an arbitrary complex
polynomial with disconnected Julia set(without the assumption that
each critical component of the filled-in Julia set is aperiodic).
We will prove that each wandering component of $K_f$ is a point.
It concludes that all but countably many components of the
filled-in Julia set are single points. This result gives an
affirmative answer of a question in a remark in Milnor's book,
see\cite{Mi}.

For any cubic polynomial with disconnected Julia set, it follows
from Theorem 5.2 and Theorem 5.3 in \cite{BH-2} that each
wandering component of the filled-in Julia set is a point.

For a polynomial with high degree, the second author proved

\begin{theoremD} [\cite{Yin}] Let $f$ be a complex polynomial
providing each critical point $c$ in wandering Julia components(if
any) is non-recurrent, i.e.,$c\notin \omega (c)$. Then each
wandering component of the filled-in Julia set $K_f$ is a point.
\end{theoremD}

For any $x\in K_f$, let
$$\mathrm{Crit}(x)=\{c\in\mathrm{Crit}\,|\,\,x\to c\},$$
where $\mathrm{Crit}$ is the set of critical points in the
filled-in Julia set.

Let
\begin{eqnarray*}
&&\textrm{Crit}_\textrm{n}(x)=\{c\in\mathrm{Crit}(x)\,| \,\, T(c) \textrm { is non-critical}\},\\
&&\textrm{Crit}_\textrm{p}(x)=\{c\in\mathrm{Crit}(x)\,| \,\, T(c) \textrm { is persistently recurrent}\},\\
&&\textrm{Crit}_\textrm{r}(x)=\{c\in\mathrm{Crit}(x)\,|\, \, T(c) \textrm { is reluctantly recurrent}\},\\
&&\textrm{Crit}_\textrm{en}(x)=\{c^\prime\in\mathrm{Crit}(x)\,|\,
\, c^\prime \not\to c^\prime
\textrm{ and } c^\prime \to c \textrm{ for some } c\in\textrm{Crit}_\textrm{n}(x) \},\\
&&\textrm{Crit}_\textrm{ep}(x)=\{c^\prime\in\mathrm{Crit}(x)\,|\,
\, c^\prime \not\to c^\prime
\textrm{ and } c^\prime \to c \textrm{ for some } c\in \textrm{Crit}_\textrm{p}(x)\},\\
&&\textrm{Crit}_\textrm{er}(x)=\{c^\prime\in\mathrm{Crit}(x)\,|\,
\, c^\prime \not\to c^\prime \textrm{ and } c^\prime \to c
\textrm{ for some } c\in \textrm{Crit}_\textrm{r}(x)\}.
\end{eqnarray*}
Then
$$\mathrm{Crit}(x)=\textrm{Crit}_\textrm{n}(x)\cup \textrm{Crit}_\textrm{p}(x)\cup \textrm{Crit}_\textrm{r}(x)\cup
\textrm{Crit}_\textrm{en}(x) \cup \textrm{Crit}_\textrm{ep}(x)\cup
\textrm{Crit}_\textrm{er}(x).$$

\begin{prop}
Suppose $x\in K_f$ and $x\not\to c$ for any critical point $c$
contained in a periodic component of the filled-in Julia set
$K_f$. Then $K_f(x)=\{x\}.$
\end{prop}
\begin{proof}
If $\mathrm{Crit}(x)= \emptyset$, then $T(x)$ is non-critical. By
the Proposition 1(1), we have $K_f(x)=\{x\}$.

If $\textrm{Crit}_\textrm{n}(x)\cup
\textrm{Crit}_\textrm{r}(x)\neq \emptyset$, by the same methods as
in the proof of Proposition 1(2), we have $K_f(x)=\{x\}$.

Now we suppose that
$$\mathrm{Crit}(x)=\textrm{Crit}_\textrm{p}(x)\cup \textrm{Crit}_\textrm{ep}(x)\neq \emptyset .$$
Since $x\not\to c$ for any critical point $c$ contained in a
periodic component of the filled-in Julia set $K_f$, hence $T(c)$
is not periodic for any $c\in \mathrm{Crit}(x)$. By the proof in
the Main Proposition, we have $K_f(x)=\{x\}$.
\end{proof}

We state a result stronger than the Main Theorem as the following

\begin{Theorem} Let $f$ be a polynomial of degree $d\geqslant 2$ with
a disconnected Julia set and let $K$ be a connected component of
the filled-in Julia set $K_f$.

(1)\; If $f^n(K)$ is a periodic component for some $n \geqslant 0$
and there is at least one critical point in the cycle of this
component, then $K$ is not a point.

(2)\; If $f^n(K)$ is a periodic component for some $n \geqslant 0$
and there is no critical points in the cycle of this component,
then $K$ is a point.

(3)\; If $K$ is a wandering component, i.e., $f^n(K)$ is not
periodic for all $n \geqslant 0$, then $K$ is a point.
\end{Theorem}
\begin{proof}
The proofs of (1) and (2) are routine, see \cite{BH-2}.

By iteration, we may assume that each periodic component
containing critical points(if any) is invariant. Let $K$ be a
wandering component of $K_f$ and $x$ be a point in $K$. Then
$K=K_f(x)=\bigcap_{k\geqslant 0}P_k(x)$. There are two
possibilities
\begin{enumerate}
\item[(a)] There is a critical point $c_0$ contained in an invariant
component of the filled-in Julia set $K_f$ such that $x\to c_0$.

\item[(b)] $x\not\to c$ for any critical point $c$
contained in an invariant component of the filled-in Julia set
$K_f$.
\end{enumerate}

In case (a), let $l_k\geqslant 1$ be  the first moment such that
$f^{l_k}(x)\in P_k(c_0)$ for any $k\geqslant 0$, i.e., $(k,l_k)$
is the first $c_0$-position on the $k$-th row in the tableau
$T(x)$. Then there is an integer $D\geqslant 1$ such that
$\deg(f^{l_k}:\; P_{k+l_k}(x)\to P_k(f^{l_k}(x)))\leqslant D$ for
all $k$. Since $K=K_f(x)$ is wandering, there exists an integer
$n_k> k$ such that $(n_k-1,l_k)$ is a $c_0$-position and
$(n_k,l_k)$ is not critical. By the tableau rule (T3) in section
2, there is no critical position on the diagonal
$$\{(n,m)\,|\,\, n+m = n_k+l_k, \quad 1\leqslant n \leqslant n_k
\}.$$ Then
\begin{align*}
\deg(f^{n_k+l_k}:\; P_{n_k+l_k}(x)\to P_0(f^{n_k+l_k}(x)))
&=\deg(f^{l_k}:\; P_{n_k+l_k}(x)\to P_{n_k}(f^{l_k}(x)))\\
&\leqslant \deg(f^{l_k}:\; P_{k+l_k}(x)\to P_k(f^{l_k}(x)))\\
&\leqslant D.
\end{align*}
There is a positive constant $\nu$ such that
$$\mod(P_{n_k+l_k}(x) - \overline{P_{n_k+l_k+1}(x)})\geqslant \nu$$
for all $k\geqslant 0$. This implies that $K=K_f(x)=\{x\}$ is a
point.

In case (b), it follows from Proposition 2 that $K=K_f(x)=\{x\}$.
\end{proof}

An immediately consequence is
\begin{corollary}
Let $f$ be a polynomial of degree $d\geqslant 2$ with a
disconnected Julia set. Then all but countably many components of
the filled-in Julia set are single points.
\end{corollary}

\begin{remark}
This corollary is not true for arbitrary rational maps, see
\cite{Mc-3}.
\end{remark}

\begin{acknowledgements} The work of this paper was deeply inspired
by a talk given by Weixiao Shen at the Morningside Center of
Mathematics in 2004. We are grateful to him for his ideas and his
patient explanations of his joint work with Kozlovski and van
Strien. We also want to express our thank to Alexander Blokh and
Guizhen Cui for their suggestions on Proposition 2 in section 5,
to Pascale Roesch and Tan Lei for their attention, suggestions and
encouragement, and to the Morningside Center of Mathematics for
its hospitality.

This work was partially supported by the National Natural Science
Foundation of China.
\end{acknowledgements}

\bibliographystyle{amsplain}

\noindent{Weiyuan Qiu}

{School of Mathematical Sciences}

{Fudan University}

{Shanghai, 200433}

{P.R.China}

{wyqiu@fudan.edu.cn}

\vskip0.2cm

\noindent{Yongcheng Yin}

{Department of Mathematics}

{Zhejiang University University}

{Hangzhou, 310027}

{P.R.China}

{yin@zju.edu.cn}

\end{document}